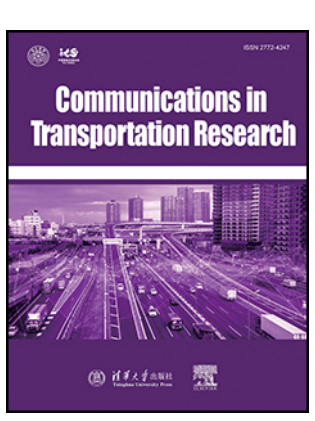

Full Length Article

# Continuum modeling of freeway traffic flows: State-of-the-art, challenges and future directions in the era of connected and automated vehicles

Saeed Mohammadian [a,b], Zuduo Zheng [b,*], Md. Mazharul Haque [a], Ashish Bhaskar [a]

[a] School of Civil and Environmental Engineering, Queensland University of Technology, Brisbane, 4000, Australia
[b] School of Civil Engineering, The University of Queensland, Brisbane, 4072, Australia

A B S T R A C T

Connected and automated vehicles (CAVs) are expected to reshape traffic flow dynamics and present new challenges and opportunities for traffic flow modeling. While numerous studies have proposed optimal modeling and control strategies for CAVs with various objectives (e.g., traffic efficiency and safety), there are uncertainties about the flow dynamics of CAVs in real-world traffic. The uncertainties are especially amplified for mixed traffic flows, consisting of CAVs and human-driven vehicles, where the implications can be significant from the continuum-modeling perspective, which aims to capture macroscopic traffic flow dynamics based on hyperbolic systems of partial differential equations. This paper aims to highlight and discuss some essential problems in continuum modeling of real-world freeway traffic flows in the era of CAVs. We first provide a select review of some existing continuum models for conventional human-driven traffic as well as the recent attempts for incorporating CAVs into the continuum-modeling framework. Wherever applicable, we provide new insights about the properties of existing models and revisit their implications for traffic flows of CAVs using recent empirical observations with CAVs and the previous discussions and debates in the literature. The paper then discusses some major problems inherent to continuum modeling of real-world (mixed) CAV traffic flows modeling by distinguishing between two major research directions: (a) modeling for explaining purposes, where making reproducible inferences about the physical aspects of macroscopic properties is of the primary interest, and (b) modeling for practical purposes, in which the focus is on the reliable predictions for operation and control. The paper proposes some potential solutions in each research direction and recommends some future research topics.

## 1. Introduction

Continuum models of traffic flow are hyperbolic systems for capturing the evolution of traffic states based only on the initial and boundary conditions (Daganzo, 1995). These models are also known as macroscopic models since they can study traffic flow's collective behaviour at an aggregated level using fluid-like state variables such as density and flow. Continuum models require less calibration effort and are computationally faster than microscopic models, and moreover, their input data are already consistent with the aggregated data from the loop detectors (Kotsialos et al., 2002; Papageorgiou, 1998).

Since the seminal contributions by Lighthill (1955) and Richards (1956), numerous developments have been made on the subject with respect to the theoretical considerations and empirical observations of real-world traffic. To be considered as desirable, models must be theoretically consistent, behaviorally sound, and capable of explaining/predicting a wide range of empirical observations, while having minimal practical and implementation issues (Helbing, 2001). However, debates and discussions are still ongoing over existing models' performances and properties, and to the best of the authors' knowledge, no existing model has yet received a strong consensus amongst researchers regarding all the aspects above.

With respect to theoretical consistency, in some cases, such as the issues around faster-than-traffic waves (Daganzo, 1995), researchers' viewpoints are diverging into different schools of thought (Helbing, 2009; Helbing and Johansson, 2009; Zhang, 2009). Meanwhile, much less has been understood about the usefulness of the proposed models against real-world traffic, since compared to the large body of literature focusing on continuum models, the empirical studies of continuum models are extremely limited. The outcomes of recent real-world evaluations of the models (Kontorinaki et al., 2017; Mohammadian et al., 2021c; Spiliopoulou et al., 2014) are conflicting regarding some aspects, such as capacity drop. Furthermore, it has been shown that some models with supposed theoretical consistency (e.g., regarding anisotropy property (Daganzo, 1995)) can exhibit significant behavioural issues of the same kind when implemented on real-world traffic (Mohammadian et al., 2021c).

* Corresponding author.
*E-mail address:* zuduo.zheng@uq.edu.au (Z. Zheng).

Issues applicable to the existing models should be revisited in the context of connected and automated vehicles (CAVs), which will share roads with conventional vehicles in the upcoming decades. Properties of such mixed traffic conditions will be much more complex because of the more pronounced role of human factors and the specific aspects of CAVs, such as information/connectivity, automation, lane choice, etc. Therefore, the underlying assumptions of the existing models for conventional traffic may no longer be valid to describe the properties of the mixed traffic condition. There is a strong need to develop new continuum models as we transition into the era of CAVs.

A focused review of continuum models is constructive in identifying the potentials and limitations of existing models as well as revealing the research needs for the development of new models. While there are already several instructive reviews on the subject (Bellomo and Dogbe, 2011; Ferrara et al., 2018 (chapters 3 and 4); Helbing, 2001; Hoogendoorn and Bovy, 2001; Seo et al., 2017; Treiber and Kesting, 2013 (chapters 6–9); Van Wageningen-Kessels et al., 2015), each one has its own specific limitations. Some notable review efforts date back to twenty years ago, and there have been substantial developments on the subject in the last two decades, whereas some recent reviews focus on theoretical aspects (Bellomo and Dogbe, 2011) or historical connections between different continuum traffic flow models (Van Wageningen-Kessels et al., 2015). More importantly, none of the existing efforts have provided a review of the recent developments in modeling CAV traffic. As such, many critical issues inherent to modeling mixed CAV traffic flows have not been explored in the context of continuum modeling.

This paper aims to highlight and discuss some essential problems in continuum modeling of real-world freeway traffic flows in the era of CAVs. The review parts of this paper put CAVs as the main ingredient of future directions. To this end, first, a select review of the major continuum-modeling families for conventional traffic is provided, focusing on those with direct implications or applications for CAV traffic flows. Along our review, we also present the models proposed for CAV traffic flows, or mixed traffic flows. Wherever applicable, we provide new insights on the properties and interpretations of the existing continuum models for traditional traffic while also discussing the implications for the CAV traffic. Our scope is restricted to vehicular interactions at the freeway level, and thus, network traffic flow theories and node models are beyond the scope of this paper (Garavello and Piccoli, 2006a; Jin, 2021).

The discussion part of this paper takes a closer look at several aspects of CAVs that are relevant to macroscopic modeling, and highlights some essential problems inherent to continuum modeling of real-world (mixed) CAV traffic flows by distinguishing between two major research directions: (a) modeling for explaining purposes, where making reproducible inferences about the physical aspects of macroscopic properties is of the primary interest, and (b) modeling for practical purposes, in which the main focus is reliable predictions for operation and control. The paper proposes some potential solutions in each research direction and recommends some future research topics.

## 2. Equilibrium models of traffic flows: applications and implications for CAVs

This section provides a select review of equilibrium continuum models developed for traffic flows of human-driven vehicles as well as some recent efforts for incorporating CAVs into continuum modeling. The phrase “equilibrium model” in this paper refers to those models where macroscopic traffic dynamics are determined, in the time–space diagram, primarily based on the mass conservation principle, where the fundamental diagram is the primary mechanism to define the relationships between macroscopic state variables (i.e., density, flow, and speed).

Our review is restricted to some select families of equilibrium models, which have some clear implications for CAV-related issues. To stay focused, we exclude some major modeling frameworks, such as LWR-type models with stochastic speed–density relationships (Bonzani and Mussone, 2002; Jabari and Liu, 2012; Kim and Zhang, 2008; Park and Abdel-Aty, 2011; Sopasakis and Katsoulakis, 2006; Sumalee et al., 2011).

We first present some basic principles and analytical and empirical findings of macroscopic traffic flow dynamics in the early decades. We then focus on some more recent continuum models in different families and elaborate on some of their behavioural and empirical properties, while revisiting the implications for (mixed) CAV traffic flows in real world.

### 2.1. Lighthill–whitham–richards model: a quick recap of some basic empirical and analytical aspects

Lighthill (1955) and Richards (1956) utilized the kinematic wave theory and proposed the first continuum traffic flow model based on the differential form of mass conservation law (LWR for short). For a homogeneous section of a road, treated as a whole, the LWR model is given as

$$\partial\rho(x,t)\,/\,\partial t+\partial\rho V(x,t)\,/\,\partial x=0 \tag{2.1}$$

where $\rho(x,t)$ is density (veh/km) and $V(x,t)=V_e(\rho)$ is the average traffic speed (km/hr) assumed to be an explicit function of density, and $Q=\rho V$ is traffic flow (veh/hr). In this case, the propagation of traffic waves can be studied by tracking the solutions to the Riemann problems (i.e., problems involving a jump in otherwise continuous initial conditions) in time and space (Lebacque (1996) for a detailed discussion of the Riemann problems).

The LWR model can replicate the basic wave-type properties of traffic flow such as congestion emergence, propagation, and dissolution, but it falls short regarding some complex empirical traffic phenomena such as: (a) the scattering phenomenon, commonly observable in the empirical flow–density diagrams (Cassidy, 1998; Nelson and Sopasakis, 1998; Schönhof and Helbing, 2009), (b) the capacity-drop phenomenon typically observable at the front of the bottlenecks, where vehicles leave the congestion zone to the free-flow condition further downstream (Hall and Agyemang-Duah, 1991; Leclercq et al., 2011; Treiber et al., 2006), and (c) the hysteresis phenomenon along the phase trajectory of the traffic states, which can be identified by tracking and connecting the macroscopic traffic states in acceleration and deceleration cycles (Treiterer and Myers, 1974).

The functional form of the fundamental diagram plays an integral part in the analytical properties of the LWR model and the corresponding physical intepretations. In order to obatin physically relevant results for traffic flow, the fundamental digram should be concave (Ansorge, 1990). If the fundamental diagram is strictly concave, the LWR model produces shockwaves in the deceleration situations and rarefaction waves in the acceleration situation (Ansorge, 1990). A strictly concave *FD* may, however, result in several unrealistic performances of the LWR model in terms of driving behaviours. For instance, a strictly concave fundamental diagram requires that vehicles always adapt their speed to the small changes in density due to the genuine nonlinearity of the characteristic speed (Zhang, 2001), and thus, the magnitude of the corresponding waves would always be variable for the entire density region. However, empirical observations have revealed that beyond certain congested states, the congestion propagation wave speed is often relatively insensitive to traffic density (e.g., Edie and Baverez, 1967; Treiber and Kesting, 2013 (chapter 8)).

### 2.2. Cell-transmission model and its extensions and implications for CAVs

Daganzo (1994) proposed the cell transmission model (CTM), which is basically a discretised cell-based variation of Eq. (2.1) coupled with a bilinear (aka triangular) fundamental diagram as

$$Q_e(\rho) = \begin{cases} \rho V_{\max}, & \text{if } \rho < \rho_{cr} \\ |C_{jam}|(\rho_{jam} - \rho), & \text{otherwise} \end{cases} \tag{2.2}$$

where $V_{\max}$ is the maximum speed, $C_{jam}$ is the propagation velocity of kinematic waves moving through the congestion region (i.e., $\rho \in [\rho_{cr}, \rho_{jam}]$) with $\rho_{cr}$ being the critical density, given from the rest of the parameters as $\rho_{cr} = \rho_{jam}|C_{jam}|/(V_{\max} + C_{jam})$ is the critical density.

From the perspective of CTM, both the kinematic waves moving through congested states and the propagation velocity of density variations in the congested states are equivalent, and of the magnitude $|C_{jam}|$. On the other hand, widespread empirical observations of conventional traffic flow's macroscopic properties have identified that beyond a certain point, the propagation velocity of kinematic waves moving through traffic flow is relatively invariant to density, which is typically in the range of [(−18)–(−12) km/hr] (Helbing et al., 2013; Treiber and Kesting, 2011, 2013; Zheng et al., 2011a; Zielke et al., 2008). As such, triangular-like fundamental diagrams have become more popular for continuum modeling (Daganzo et al., 1997; Del Castillo, 2001, 2012; Del Castillo and Benítez, 1995b).

The triangular fundamental diagram is the continuum approximation of the steady-state solution of Newell's car-following model (Newell, 2002) in which case, the average time gap can be obtained as $\tau_{gap} = (\rho_{jam}|C_{jam}|)^{-1}$, which is constant across $[\rho_{cr}, \rho_{jam}]$. This property seems appealing for describing the steady-state traffic flow of CAVs. For instance, the constant time gap is equivalent to a constant spacing policy. Some recent empirical studies of CAV traffic have shown the existence of a triangular fundamental diagram in the steady-state traffic using empirical trajectories collected from field testing with CAVs (Li et al., 2022b; Shi and Li, 2021). These studies have also identified and formulated time-gap specific fundamental diagrams for commercial CAVs' steady-state traffic flow, where the maximum capacity corresponds to the shortest headway setting (Li et al., 2022b; Shi and Li, 2021).

Meanwhile, several recent efforts have been made to extend CTM for continuum modeling of traffic flow of CAVs (Jin and Yang, 2013; Jin et al., 2022; Wu et al., 2022; Zhou and Zhu, 2020). Such studies' outcomes indicate that as CAV market penetration rises or time gaps shorten, traffic capacity increases and congestion regions are alleviated. These findings align with empirical studies (Li et al., 2022b; Shi and Li, 2021) which show vehicles equipped with adaptive cruise control (ACC) having greater capacity at smaller time-gaps.

Some issues must be considered regarding using CAVs' fundamental diagrams in the steady-state traffic to make inferences about the kinematic waves in CAV traffic flows. For instance, considering the microscopic interpretation of the cell-transmission model, i.e., the model by Newell (2002), a constant time gap in the triangular fundamental diagram is considered to be the same as the drivers' response time, i.e., the adaptation time to changes in the traffic condition ahead (stimulus) at all density regimes, regardless of the nature of the stimulus (e.g., acceleration or deceleration) and its intensity (speed difference, etc.). While such an assumption is inadequate even for human-driven traffic (Zhang, 1999), its inadequacy for CAVs is likely to be much more pronounced for a variety of reasons, including the fact that safety considerations and proactive rear-end crash avoidance are explicitly incorporated in the CAVs' underlying control laws (Calvert et al., 2018; Delle Monache et al., 2019; Shladover, 2018). Such considerations largely depend not only on the nature of the stimulus ahead, but also on its intensity regarding several factors such as CAVs' speed, speed difference, time headway, time-to-collision, etc. (Hicks, 2018; Shladover, 2005). Section 4.1.2 provides a more comprehensive discussion on this subject using empirical and theoretical findings in the literature.

### 2.3. LWR model with bounded acceleration

One of the main drawbacks of the original LWR model is its unrealistic speed adaptation mechanism. To elaborate, the macroscopic speed adaption of the vehicles is considered, which is

$$\mathrm{d}V_e(t)/\mathrm{d}t = -\rho\left(V_e'(\rho)\right)^2 \partial\rho/\partial x \tag{2.3}$$

Equation (2.3) imposes shocks near discontinuities, suggesting that speed adaptation rate is unbounded both for deceleration and acceleration maneuvers (if the fundamental diagram is triangular-like (Zhang, 2001)).

The first effort to incorporate bounded acceleration into the LWR model was made by Lebacque (1997). The author proposed some adjustments for the numerical solution of the LWR model, where a systematic bound is applied for the cells experiencing unbounded acceleration. In Lebacque (2002), a two-phase extension of the LWR model with bounded acceleration (BA-LWR for short) was developed, which directly incorporated boundedness of the acceleration into the governing equations, by considering an upper bound. If the acceleration predicted by the LWR model is beyond the maximum, a dynamic equation for speed ($\mathrm{d}V/\mathrm{d}t = A_{\max}$) is coupled with the LWR model, where $A_{\max}$ is the bounded acceleration rate. The BA-LWR model can be presented as

$$\begin{cases} \partial\rho\Big/\partial t + \partial\rho V\Big/\partial x = 0, V = V_e(\rho), \quad \text{if } \dfrac{\mathrm{d}V_e(\rho)}{\mathrm{d}t} \le A_{\max}, \text{ LWR phase} \\ \begin{cases} \partial\rho/\partial t + \partial\rho V/\partial x = 0, \\ \partial V/\partial t + V\partial V/\partial x = A_{\max}, \end{cases} \text{otherwiswe, BA phase} \end{cases} \tag{2.4}$$

Note that the BA-LWR model is technically non-equilibrium if the BA phase is activated. However, in this paper we classify the BA-LWR model as equilibrium because the BA phase activates only in certain situations.

It has been shown that with a triangular fundamental diagram, BA-LWR is equivalent to an extension of the car-following model by Newell (2002) for the bounded acceleration (Jin and Laval, 2018). The BA approach has been utilized to capture various aspects of traffic flow within the LWR framework. Leclercq (2007) utilized the BA-LWR framework for studying moving bottlenecks and argued that the effects of slow-moving vehicles were underestimated in previous studies of moving bottlenecks within the LWR model (e.g., Newell, 1998) due to the unbounded acceleration. It has been argued that the BA-LWR model by Leclercq (2007) can also be applied for noise and emission modeling of moving bottlenecks (Can et al., 2010). The approach proposed by Leclercq (2007) applies the bounded acceleration to almost all vehicles in a single-class traffic flow, whereas moving bottlenecks are often caused by one or a group of slow-moving vehicles. As well, slow moving and heavy vehicles are the major causes of factors such as noise and emission. Therefore, we argue that a bounded acceleration framework applied to single-class equilibrium models may not have significant implications for factors such as moving bottlenecks and noise and emission modeling.

Laval and Daganzo (2006) utilized the BA approach to accommodate lane-changing traffic flows on multi-lane freeways within the BA-LWR framework. In this work, the LWR model is applied for each lane separately, where lane-changing vehicles are defined as moving bottlenecks subject to bounded acceleration. Through numerical simulations, the authors argued that the proposed model could capture the capacity-drop phenomenon observed downstream of lane-drop bottlenecks and be useful for studying the capacity flow of moving bottlenecks. Jin (2010) argued that lane-changing vehicles hinder traffic flow not only in the destination lane but also in their current lane. To account for such effects,

Jin (2010) proposed another LWR-type model in which the effect of lane-changing vehicles in all lanes is reflected as a capacity drop in the fundamental diagram. The model was validated in Gan and Jin (2013). Jin (2017a) developed a model for lane-drop bottlenecks and analytically investigated the properties of kinematic waves in lane-drop scenarios. The model was further developed in Jin (2017b) to account for capacity drop, where lane-changing vehicles inside lane-drop zones are assumed to adopt a bounded acceleration rate.

It has been argued that the BA-LWR model can capture capacity drop and hysteresis (Khoshyaran and Lebacque, 2015). While numerical simulation studies have shown that the BA-LWR framework can replicate the capacity drop phenomenon, controversies remain around the performance of the BA-LWR model with respect to capacity drop, when implemented on real-world traffic through calibration. Our recent benchmarking study of the continuum model against real-world traffic (Mohammadian et al., 2021c) shows that the BA-LWR model may have a limited capability of capturing the observed capacity drop and hysteresis due to several reasons. For instance, the BA-LWR model only captures the capacity drop phenomenon through the difference between acceleration and dismisses other potential contributing factors, such as the lane-changing maneuvers in the merge zones or the proportion flow of merging/diverging vehicles. Therefore, since in the BA-LWR model the magnitude of capacity drop is a property of the upper acceleration bound, realistic values for the acceleration bound may cause unrealistic values for capacity drop and vice versa.

Now let us have a closer look at the implications of these points and the BA-LWR model for the traffic flow of CAVs. One may pick up the BA-LWR framework to investigate to what extent the presence of CAVs in traffic flow may affect the capacity drop phenomenon. This problem is similar to the investigation performed by Jin and Laval (2018) for the human-driven traffic, in which a unified BA-LWR framework was proposed. The authors investigated the impact of variable speed limit on the capacity drop phenomenon, using their proposed BA-LWR framework and concluded that bounded acceleration does not lead to capacity drop within the variable speed limit zones. Since CAVs can be used as moving controllers to alleviate undesirable traffic flow properties, one may adopt an approach similar to the one proposed by Jin and Laval (2018) in order to implement the variable speed limit through CAVs, in which case, the findings may suggest that CAVs can remove or significantly reduce the capacity drop phenomenon when approaching the bottleneck zones from behind. However, such findings are likely to be inclusive due to the inherent limitations of the BA-LWR framework in capturing the capacity drop phenomenon comprehensively, as discussed above.

### 2.4. Issue of bounded deceleration and the LWR-type with non-local speed-density relationship

While bounded speed adaptation rate is an issue for both deceleration and acceleration maneuvers, most studies have focused on incorporating the boundedness of the acceleration rate into the continuum models (e.g., Jin and Laval, 2018; Laurent-Brouty et al., 2021; Lebacque, 2002; Lebacque, 1997; Leclercq, 2007a), possibly for two reasons. First, for human-driven traffic, the boundedness of the acceleration is a more important issue; and the magnitude of acceleration rate is typically restricted to a small range ($\sim 2\ \mathrm{m/s^{-2}}$ (Treiber and Kesting, 2013)), whereas sharp deceleration rates ($> 8\ \mathrm{m/s^{-2}}$ (Sharma et al., 2019a) have been identified in empirical observations of human-driven vehicle trajectories.

However, perhaps the more important reason is the inherent difficulties that incorporating bounded deceleration might bring to the mathematical tractability of the continuum models (Giorgi et al., 2002). Several consistency conditions for hyperbolic systems, such as the Lax inequalities and the Rankine–Hugoniot jump condition (LeVeque, 2002), may not be satisfied, if the deceleration rate near discontinuities is bounded.

When it comes to the traffic flow of CAVs, boundedness of deceleration could become a more significant and more important ingredient of such traffic. Empirical observations of existing field data with CAVs have identified CAVs generally drive more safely and have a smoother deceleration rate (Hu et al., 2023). To the best of the authors' knowledge, the issue of shocks and unbounded deceleration has not endogenously addressed in the classical LWR-type models.

Meanwhile, another non-classical LWR-type model, namely LWR-type models with non-local speed-density relationships (Chiarello, 2021; Chiarello and Goatin, 2018a, 2018b, 2019; Goatin and Scialanga, 2015, 2016; Karafyllis et al., 2022; Kurganov and Polizzi, 2009; Lee and Liu, 2013; Li and Li, 2011; Sopasakis and Katsoulakis, 2006), in which drivers' preferred speed depends, not on the local density, but on a range of traffic condition in the visible range. Some of these models can be expressed in a generic form as

$$\partial \rho \,/\, \partial t + \left[ \partial \rho V \left( \int_x^{x+\eta} \rho(y,t) \omega_\eta (y-x) \mathrm{d}y \right) \right] \Big/ \partial x = 0 \tag{2.5}$$

where $\eta$ is a constant and represents drivers' visible range, $\omega_\eta$ is a kernel that determines drivers' speed at $x$, based on a weighted average of traffic conditions in the range of $[x, x+\eta]$. The well-posedness of these non-local models have been proven (Keimer and Pflug, 2017). These models seem to have promising potentials such as eradicating shocks, and smooth and asymmetric speed adaptation rates, and therefore, have been recently utilized for modeling traffic flow of CAVs (Chiarello, 2019; Guan, 2022; Huang, 2022; Huang and Du, 2022). However, several behavioral and physical considerations need to be taken into account and thoroughly investigated before making inferences about CAV traffic flows using non-local LWR-type continuum models:

(1) The analytical and numerical experiments with such models show that in many cases, the initial data (either discontinuous or Lipschitz continuous) can still eventually lead to the formation of shocks (e.g., Chiarello, 2021; Lee and Liu, 2013), and thereby, unbounded deceleration. Examples of such conditions have been analytically proven by Keimer and Pflug (2017) and numerically shown in other studies (Friedrich et al., 2018; Keimer and Pflug, 2017). Therefore, we argue that any continuum model that leads to the formation of shockwaves cannot endogenously capture the bounded deceleration property.

(2) The potentials for asymmetric speed adaptation behavior such that acceleration rate becomes higher than deceleration rate, which is behaviorally unrealistic, and can result in counter-clockwise hysteresis loops, which while occasionally observable, are counteractive when produced systematically (Mohammadian et al., 2021c).

(3) The possibility of systematic back-to-front acceleration for a line of vehicles such that some preceding vehicles react to the changes of traffic condition further ahead more intensely than their immediate leader would, and such a scenario is unrealistic in traffic flow (Zhang, 2009).

(4) Emergence of bumps and instabilities in the density profiles of an otherwise initially homogeneous free-flow condition if the upstream traffic is heavily congested. If that is the case, it would be a strong conflict with world-wide observations of jam front dissolution. Empirical observations consistently show that vehicles leaving congestions to very-light traffic ahead quickly recover their free-flow speed, and density profile shows little dispersion (Treiber and Kesting, 2013 (chapter 8)).

(5) Some recent works demonstrate that in non-local traffic flow models, corresponding fundamental diagram may lose the concavity condition (Sun and Tan, 2020), which in turn, may produce many unrealistic behavioral properties (like deceleration fans) (Zhang, 2009).

### 2.5. LWR-type models involving multiple user/vehicle classes

The original LWR model assumes a single user class, neglecting vehicle and user heterogeneities and their effects on traffic flow. Over time, efforts have been made to include these heterogeneities, leading to multi-class LWRs (MC-LWRs). As stated by Logghe and Immers (2003), a general formulation of MC-LWR models is

$$\begin{cases} \partial\rho_u/\partial t + \partial(\rho_u V_u)/\partial x = 0 \\ V_u = V(\rho_{1,\dots,}\rho_{nu}) \end{cases} \tag{2.6}$$

In multi-class models, each vehicle class, identified by subscript $u$, has its own speed. Equation (2.6) shows MC-LWR models delineating the LWR model for each class, with classes interlinked through the fundamental diagram. The main variation in these models is in defining these relationships.

Some MC-LWR models account for distinct route-choice behaviors per class and can be termed multi-commodity models, suitable for traffic networks (e.g., Jin, 2012, 2013). Daganzo et al. (1997) introduced a two-user-class MC-LWR model, to accomdate the queuing effects, in which one class has a designated lane. This model mimics the behavior of vehicles diverging near off-ramps and their impact on traffic.

Daganzo (2002) introduced an MC-LWR model distinguishing between two user classes on a two-lane highway: "rabbits" (aggressive drivers) and "slugs" (timid drivers). Rabbits overtake often, using high-speed lanes, while slugs avoid passing. The model assumes drivers' behaviors remain consistent regardless of traffic conditions. Two traffic regimes are defined: "1-pipe" and "2-pipe". The "1-pipe" regime treats multi-lanes as a single lane when traffic in both lanes moves at identical speeds. The "2-pipe" state sees varied conditions in lanes, with one class in free-flow and the other in congestion. Daganzo (2002) stated that this model could capture the reversed lambda pattern seen in flow-density data.

Wong and Wong (2002) introduced an MC-LWR model where user classes, though of identical vehicle length, had differing desired speeds. This allowed classes to adopt different speeds based on spacing. The model was believed to mimic capacity-drop, hysteresis, and platoon dispersion. Zhang and Jin (2002) presented an MC-LWR model differing in vehicle lengths and desired speeds, with varied speeds only in free-flow conditions; all vehicles move at the same speed when congested. Chanut and Buisson (2003) used similar logic in their model, with dynamic maximum and critical densities based on traffic composition. Some MC-LWR models base class dynamics on passenger cars. Logghe and Immers (2003) introduced a model scaling the effects of heavy vehicles using fixed passenger-car equivalents (PCEs), allowing varied spacing for vehicle classes at identical speeds. PCEs remain constant across traffic regimes. Ngoduy and Liu (2007) also employed PCEs in their model, where classes have different free-flow speeds but synchronize in congestion. This model is believed to capture capacity drops, hysteresis, and platoon dispersion. However, an important question is whether discontinuities in flow-density diagrams from MC-LWR models represent a capacity drop. Researchers attribute such discontinuities to interactions between vehicle classes (e.g., Ngoduy and Liu, 2007; Wong and Wong, 2002).

Van Lint et al. (2008) introduced the FASTLANE MC-LWR model, utilizing the PCE concept. In this model, PCEs, determined by speed, spacing, and temporal headway of classes, define the effective density. This effective density then dictates the speed for each class. In FASTLANE, both PCEs and the fundamental diagram are interdependent. Van Wageningen-Kessels et al. (2014) proposed a generic version of FASTLANE, which can be adapted to earlier MC-LWR models based on parameter choices. However, discrepancies in flow-density diagrams are not due to non-equilibrium traffic conditions, as all states can align with equilibrium curves based on traffic composition.

Overall, it has been argued, based on numerical simulations, that the MC-LWR type models like the above can capture some of the complex traffic phenomena such as scattering, capacity drop, and hysteresis (Ngoduy and Liu, 2007; Wong and Wong, 2002). We shall take a closer look at these aspects by contrasting the numerical results obtained by the aforementioned studies to the models' underlying equations, and real-world implementations. Regarding scattering, our recent benchmarking studies of continuum models against real-world traffic suggests that multi-class can only account for a tiny range of the scattering observed in real data (Mohammadian et al., 2021c). This observation could be explained by the fact that in MC-LWR models, simulated traffic states in the flow-density diagram can only be placed on two or a few curves of the fundamental diagrams (FDs).

Meanwhile, regarding the capacity drop phenomenon, although there is no consensus regarding the factors affecting capacity drop, the reason why the LWR model cannot describe it is an established fact due to the equilibrium traffic and symmetric phase trajectories over the FD. In other words, traffic flow for a given transition from congestion to a free-flow zone is at the maximum, which is the same as the transition from a free-flow to a congestion zone. MC-LWR models also apply the entropy condition for vehicles within each class. However, in the majority of these models, the speed of vehicle classes is different only in free-flow traffic. Therefore, when vehicle classes are considered altogether, the average flow upstream of a congestion zone could be slightly higher compared to downstream because faster vehicles reach congestion upstream sooner than slower ones. Such a flow drop could stem from averaging, which is distinctively different from capacity drop frequently reported in the literature, which can happen in traffic flow with a single-vehicle class because of the difference between acceleration and deceleration rates.

Similarly, MC-LWR models may not be suitable for studying the hysteresis phenomenon in traffic flows. It has been argued that the MC-LWR models show some closed loops in the phase diagram of the upstream of the bottleneck (Van Wageningen-Kessels et al., 2015). However, the loops produced by the MC-LWR models are not necessarily the genuine hysteresis loops because they cannot be explained by the difference between acceleration/deceleration rates. This is because MC-LWR models all belong to the equilibrium framework and their speed adaption mechanism depends primarily on the continuity equation. Furthermore, the phase trajectories predicted by MC-LWR models are collected from fixed-point locations, which means that variable groups of vehicles with different compositions of vehicle classes are compared with one another. It is worth noting that even if such loops are observable on the trajectory of moving observers; their patterns may not necessarily be in line with observed real-world data due to the role of single-pipe treatment in changing the structure of lane-based hysteresis loops.

Finally, it has been widely acknowledged that these complex phenomena are consequences of traffic instabilities in non-equilibrium regimes (Newell, 1962). Recent studies have suggested that the role of human factors (e.g., aggressive/timid driving and task difficulty) in the emergence of these phenomena is significant (Laval, 2011; Laval and Leclercq, 2010; Saifuzzaman et al., 2017; Tampère et al., 2005). MC-LWR models cannot account for either non-equilibrium traffic conditions or the impact of human factors, as the heterogeneity in these models is mainly defined as the operational properties of vehicle classes such as the desired speed and the gap acceptance behavior, and all classes have the same acceleration behavior as the standard LWR model.

The above-mentioned issues, although not directly relevant to CAVs, can become considerable when the focus is to utilize MC-LWR models to investigate the extent to which such complex traffic phenomena could be affected by CAVs.

Meanwhile, Logghe and Immers (2003) adopted a distinct approach and derived their MC-LWR model from the user-equilibrium theory (Wardrop, 1952). Their model is derived from two major premises, i.e., that vehicles do not occupy more than the necessary space and those faster vehicles do not affect the slower ones. Because of the latter, the semi-congested regime appears, in which slower vehicles can still drive at their free-flow speed, whereas faster vehicles are already in the congestion state. This model was recently further developed and generalized by Qian et al. (2017) for more than two vehicle classes, and the lateral interactions of vehicles were incorporated as the concept of perceived equivalent density. The model uses lateral interactions to determine how each class perceives the effects of other classes for a given traffic condition. It has been argued that this approach can also be useful for modeling moving bottlenecks (Logghe and Immers, 2008). However, the existence of class-specific fundamental diagrams might be open to debate, as a distinct observation of steady-state traffic for each class may not be easily obtained from detector data. Even if complete vehicle trajectories are available, the percentage of heavy vehicles experiencing steady traffic states, for a wide range of traffic densities, could be too limited to postulate a valid fundamental diagram.

Meanwhile, Qin et al. (2021) proposed a MC-LWR type model to incorporate mixed traffic flow of human-driven and CAVs, where the properties of each class were extracted from car-following relationships, such as the intelligent driver model (Treiber et al., 2000). In their proposed model, the fundamental diagram is the most conservative for human-driven vehicles, and assigns larger spacings for a given speed, whereas CAVs have the least conservative fundamental diagram, i.e., the smallest spacing for a given spacing. Using numerous hypothetical experiments, the authors investigated the effects of CAVs on traffic flows and argued that CAVs can lead to congestion alleviation and improved traffic capacity. By comparing the numerical experiments provided by the authors, it can be concluded that CAVs generally have the shortest response time, and correspondingly, the largest wave speed when leaving congestion zones, and thereby, contributing to faster congestion dissolution. As discussed previously, however, such assumptions are in sharp contrast with recent naturalistic field observations with CAVs (Hu et al., 2022, 2023; Sun et al., 2020b), revealing that due to safety considerations, CAVs' driving strategies are even more conservative than that of an average human-driven vehicle.

Other considerations need to be taken into account as well when making inferences about CAV traffic flows using existing MC-LWR type models. For instance, the majority of existing MC-LWR models are single-pipe and treat the multi-lane traffic as if it were single-lane. Consequently, in all single-pipe multi-class models, faster vehicles either can overtake the slower ones or queue behind the slowest class. In the former case, the violation of the first-in-first-out property occurs, and the impact of the lane-changing maneuvers is not adequately modelled. In the latter case, the advantage is to allow for modeling moving bottleneck. However, this property will be inadequate for modeling freeways with lane restrictions where moving bottlenecks caused by slow vehicles are restricted to specific lanes, whereas single-pipe treatment causes all the lanes to experience the moving bottleneck.

Additionally, it is likely that at early stages, road authorities will allow the drivers of CAVs to enable automated driving mode only on specifically allocated lanes (Mahmassani, 2016), such that in case of driver take over for manual driving, exiting such allocated lanes will be mandatory. On this note, some recent analytical investigations of mixed CAV traffic flow have proven that under certain conditions, lane-based traffic flow can show counter-intuitive behavior, in the sense slower vehicles tend to remain on slower lanes (Li et al., 2022a). Under such scenarios, the assumption of single-pipe traffic will no longer be adequate due to varying traffic conditions and different underlying premises across the lanes. As such, the benefits of the single-pipe models (e.g., parsimony and simplicity) may no longer outweigh the complexities of the multi-lane framework, such as designing the lane-changing mechanisms.

Therefore, it is expected that multi-lane models, instead of single-pipe models, are likely to attract more research attention for modeling the mixed traffic flow of CAVs, especially given that some of the complexities associated with the multi-lane framework might be less challenging than before. For instance, one of the main challenges used to be a lack of reliable data and adequate knowledge about the mechanisms behind discretionary lane-changing maneuvers. In recent years, on the other hand, there has been significant development in understanding lane-changing behaviors (e.g., Ali et al., 2019; Ali et al., 2020; Zheng, 2014; Zheng et al., 2013) with a focus on human factors. These understandings could be utilized for designing more reliable lane-changing mechanisms in the continuum framework.

## 3. Non-equilibrium models of traffic flow: Applications and implications for CAVs

Non-equilibrium continuum models (aka second-order models) are those in which one or multiple hyperbolic partial differential equations (PDEs) are coupled with the continuity equations (Eq. (2.1)). In non-equilibrium models, the fundamental diagram only applies to the equilibrium states in the steady-state traffic given $\partial\rho/\partial t \approx 0$, $\partial\rho/\partial x \approx 0$, $\partial V/\partial t \approx 0$, and $\partial V/\partial x \approx 0$.

There are three conventional approaches for developing non-equilibrium models:

(1) The phenomenological approach in which hyperbolic conservation laws are directly borrowed from fluid dynamic concepts and implemented for traffic flow. In this case, the speed adaptation mechanisms may not necessarily be meaningful with respect to driving behaviors (Treiber and Kesting, 2013).
(2) Derivation from car-following relations using certain approximation of car-following variables in the Eulerian coordinates (Helbing and Johansson, 2009). In this case, it may or may not be possible to establish a direct relationship between the underlying car-following models and the corresponding non-equilibrium (Papageorgiou, 1998). For instance, the car-following model by (Newell, 2002) is equivalent to the LWR model with triangular fundamental diagram (i.e., CTM (Daganzo, 1994; Daganzo, 2006a). However, with alternative approximations of car-following variables, one can derive the non-equilibrium model by Payne (1971), from the car-following model by (Newell, 2002).
(3) Derivation from the gas-kinetic approach, where vehicular interactions for a finite number of vehicles at the mesoscopic scale are transformed to the macroscopic scale (Helbing, 1996b). In this case, the speed adaptation mechanism does not capture the dynamics of any individual vehicle, but it represents the average vehicular interactions within a moving window travelling with the average traffic speed (Treiber and Kesting, 2013). The gas-kinetic theory has also been utilized to provide physical foundations to some existing car-following-based and phenomenological models (Ben-Naim and Krapivsky, 2003; Helbing, 1996b; Marques Jr and Méndez, 2013; Phillips, 1979). Helbing (1996b) showed that a hierarchy of continuum models, including the LWR and existing PW-like models at the times, could be consistently derived from the model by Paveri-Fontana (1975), which is based on the gas kinetic theory (GKT).

In the following, we provide a select review of continuum models.

### 3.1. Non-equilibrium model developed for conventional human-driven traffic

#### 3.1.1. Elementary non-equilibrium models and Daganzo's critique

The first non-equilibrium continuum model was proposed by Payne (1971), and later on independently by Whitham (1974) (known as the PW model for short). The PW model defines acceleration as a function of both deviations from the equilibrium state and changes in the traffic condition ahead. The acceleration mechanism in the PW model is

$$\partial V / \partial t + V\partial V / \partial x = (V_e(\rho) - V) / \tau - C^2(\rho) / \rho(\partial\rho / \partial x) \tag{3.1}$$

where $C(\rho) = \sqrt{-V_e'(\rho)/2\tau}$ is the so-called sonic velocity, and $\tau$ is the relaxation time.

The PW model allows for the occurrence of a wide range of traffic flows at a given density and vice versa. As a result, the scattering in the flow–density diagram can be simulated (Payne, 1979). It has been argued that the model can also replicate the capacity-drop phenomenon since the assumption of equilibrium traffic is relaxed (Papageorgiou, 1998). It has been argued that the model can also replicate the capacity-drop phenomenon (Papageorgiou, 1998).

It has been argued that the PW model is numerically challenging to simulate (Treiber and Kesting, 2013 (chapter 9)). Inconsistent accounts of the performance of the model for similar scenarios (e.g., Cremer and May 1987; Leo and Pretty, 1992; Papageorgiou, 1983) have also been linked with the choice of numerical schemes (Zhang, 1998). To overcome these issues, Messmer and Papageorgiou (1990) proposed METANET, which is based on a simplified discretised form of the PW model. METANET has been widely used for control purposes because of its simplicity (e.g., Ferrara et al., 2018 (chapter 8), and references therein). It is worth noting that although METANET is commonly considered as a discretised version of the PW model, it is a phenomenological model, and cannot be drived from car-following relations.

Depending on the choice of relaxation time, solutions of the PW model for initially homogeneous traffic subject to small perturbations can be unconditionally stable, which means that the model cannot adequately replicate traffic instabilities and stop-and-go waves (Helbing, 1996a). To capture traffic instabilities, Kühne (1984) and Kerner and Konhäuser (1993) modified the PW model and assumed that the sonic velocity in the anticipation term is constant. They also introduced a concave–convex fundamental diagram to replicate stop-and-go waves in a specific density range of congested traffic. Consequently, the stability of solutions is guaranteed not only for light traffic but also for highly congested traffic above a specific density value. The acceleration mechanism in this model is

$$\partial V / \partial t + V \partial V / \partial x = (V_e(\rho) - V) / \tau - {C_0}^2 / \rho(\partial \rho / \partial x) + \eta\left(\partial^2 V / \partial x^2\right) \quad (3.2)$$

where $C_0$ represents the so-called sonic velocity, and $\eta$ is the viscosity coefficient. The third term on the right-hand side of Eq. (3.2) is called the viscosity (or diffusion) term and is derived from a translation of Navier–Stokes' fluid dynamics equations (Helbing, 1996b) for more discussion on Navier-Stokes' equations). It has been argued that viscosity term is not physically meaningful for traffic flow but back then it would be used mainly to improve the numerical properties of the model and smooth the shocks (Treiber and Kesting, 2013 (chapter 9)). It is worth noting that today, the numerical solution of continuum models is less of a challenge as improved numerical schemes have been presented (Delis et al., 2014; Helbing and Treiber, 1999; Jin and Zhang, 2001; Mammar et al., 2009; Mohammadian and Van Wageningen-Kessels, 2018).

In preceding decades, several attempts have been made to investigate some properties of CAV-related aspects such as adaptive cruise control using PW-type models or their extensions. For instance, Darbha and Rajagopal (1999) investigated the effects of vehicle control laws on traffic flow stability by accommodating the constant time headway policy into the PW model (Payne, 1971; Whitham, 1974) and concluded that such a policy will always lead to unstable traffic flow. Later on, Li and Shrivastava (2002) argued that the findings in Darbha and Rajagopal (1999) could be biased in that the external factors in the numerical experiments, such as flow exchanges at the boundaries, were not controlled. To eliminate this effect, they investigated traffic flow behavior on a homogeneous circular road, and the numerical experiments and stability analysis suggested improved stability in the traffic flow. Another study by Yi and Horowitz (2006) derived general stability criteria for traffic flow in the presence of Adaptive Cruise Control (ACC)-equipped vehicles through performing a non-linear stability analysis on the extended version of the PW model.

Daganzo (1995) called into question the consistency of non-equilibrium models at the time with traffic flow and driving rules in several ways, including:

- **Backward driving or the wrong-way-travel effect:** The central diffusion terms added to smooth shock waves (as in the model by Kerner and Konhäuser (1993)) leads to dispersion of vehicles from behind if traffic is lighter in the upstream direction.
- **Violation of the anisotropy property:** The second characteristic speed in these models travels faster than the average traffic speed, carrying the information from the upstream traffic to the downstream vehicles, causing them to slow down or speed up accordingly, depending on the traffic condition behind, and such a property is counter-intuitive with respect to car following.

#### 3.1.2. *Some improved non-equilibrium models with respect to Daganzo's critique*

Several efforts have been made to address these criticisms by Daganzo (1995). In this section, we review some of these works, where the implications are mainly for conventional traffic.

Zhang (1998) showed that the anticipation mechanism of the PW model can itself produce the "wrong-way travel" effect, regardless of the presence of the viscosity term. The observation by Zhang (1998) is consistent with the zero-relaxation time limit of the PW, which results in the LWR model with a diffusion term (Treiber and Kesting, 2013). To resolve the issue, Zhang (1998) derived the PW model from an improved car-following relationship in which drivers respond to the frontal stimuli with a delay. The acceleration mechanism in the resulting non-equilibrium model is

$$\partial V / \partial t + V \partial V / \partial x = (V_e(\rho) - V) / \tau - C(\rho)^2 / \rho(\partial \rho / \partial x) \quad (3.3)$$

where $C(\rho) = \alpha\rho V_e^{'}(\rho)$ being the generalized form of the sonic velocity, with $\alpha$ as a constant, is the main difference between this model and the standard PW model. Zhang (1998) showed that such a modification prevents the "wrong-way-travel" effect.

Aw and Rascle (2000) argued that faster-than-traffic waves appeared in previous models because the anticipation term was derived from the standpoint of a fixed-point observer and did not take the convection of vehicles into account. They suggested that the anticipation term must be computed from the standpoint of the moving observer travelling with the average traffic speed. Such a rationale results in the following acceleration mechanism:

$$\partial V / \partial t + V \partial V / \partial x = \rho P^{'}(\rho) \partial V / \partial x \quad (3.4)$$

where $P(\rho)$ is a smooth increasing function, known as the so-called pressure variable. The characteristics of this model are $\lambda_1 = V - \rho P^{'}(\rho)$ and $\lambda_2 = V$, where the second characteristic always travels at the speed of the traffic. Therefore, the model is anisotropic in Daganzo's (1995) sense. The original Aw-Rascle model is phenomenological and is not derived from car-following relationships. However, Aw et al. (2002) showed that the Aw-Rascle model can be viewed as the hydrodynamic limit of "follow-the-leader'' (also known as stimulus-response) car-following models (e.g., Gazis et al., 1961).

Zhang (2002) independently derived a similar continuum model from a car-following relationship, assuming that drivers' response time is an explicit function of the inter-vehicular spacing. A direct consequence of this assumption is that the relaxation term does not appear in the resulting continuum model, which is in line with the model by Aw and Rascle (2000). It can be shown that Zhang's model is a special case of Aw and Rascle's (2000) model with the condition $P(\rho) = -V_e(\rho)$. These models are collectively referred to as the Aw-Rascle-Zhang (ARZ) model. As the relaxation term is absent in the ARZ model, the fundamental diagram cannot be uniquely derived since according to Eq. (3.2), any homogeneous state satisfying $\partial V / \partial x = 0$ could be viewed as an equilibrium steady-state condition. Lebacque et al. (2007b) argued that the ARZ model can be viewed as the LWR model with "variable fundamental diagram".

Although the model by Zhang (2002) incorporates the concept of the fundamental diagram, the ARZ model does not explicitly consider a maximum density, and the model can simulate an initial condition where vehicles are bumper-to-bumper but can travel at any speed. Greenberg (2002) coupled the ARZ model with a relaxation term to resolve this issue so that the non-physical solutions can approach the fundamental diagram. However, this approach addresses the issue only partially and requires a small relaxation time to be effective.

Lebacque et al. (2007b) proposed the generic "second-order" model (GSOM), which includes both the ARZ and Colombo's (2002) model as special cases. GSOM is also a phenomenological model which can include an arbitrary source term as a special case. The second dynamic equation in GSOM in non-standard form as

$$\partial(\rho\mathscr{L}) / \partial t + \partial(\rho\mathscr{L}V) / \partial x = \varphi f(\mathscr{L}) \tag{3.5}$$

where $\mathscr{L} = \mathscr{L}(\rho, V)$ is an arbitrary conserved variable satisfying the condition $\mathrm{d}\mathscr{L}/\mathrm{d}t = 0$, meaning it is conserved along its trajectory. The parameter $\varphi \in \{0,1\}$ allows one to use the model with or without the source term. It should be mentioned that the source term in this model does not, in general, equal the relaxation term as in previous PW-type models. This is because the generic quantity $f(\mathscr{L})$ does not necessarily represent the relaxation acceleration. However, with $\mathscr{L} = (V - V_{\mathrm{e}}(\rho))$, and $f(\mathscr{L}) = -\rho\mathscr{L}/\tau$, and $\varphi = 1$, GSOM is identical to the model by Greenberg (2002) and the model by Zhang (2002) with a relaxation term. Likewise, the model is the same as Colombo's (2002) model with $\mathscr{L} = (q - q_*)/\rho$ and $\varphi = 0$. Lebacque and Khoshyaran (2013) presented the Lagrangian formulation of GSOM from the perspective of a moving observer. Zhang et al. (2009) proposed another GSOM-like model in which the second equation is based on an artificial variable, called "pseudo density", which is transformed from speed.

The ARZ model has been interpreted as the LWR model with a "variable fundamental diagram", and several studies have suggested that such a GSOM model can theoretically explain the wide scattering of the traffic states by placing them through multiple flow-density curves (Fan and Seibold, 2013; Lebacque et al., 2007b; Würth et al., 2022). However, this theoretical property may not be achieved when the model parameters are estimated from minimizing the error between simulation and observation. Our recent benchmarking studies of continuum models against real-world traffic show that GSOM models can have limited performance in replicating the observed scattering (Mohammadian et al., 2021c).

This observation could, in part, be related to the model's stability behavior. A simple linear stability analysis following the procedure in Helbing and Johansson (2009) shows that the most popular variant of GSOM models, i.e., the model by Zhang (2002) with a relaxation term, is unconditionally stable and cannot replicate the emergence of stop-and-go waves from small disturbances in otherwise steady-state equilibrium traffic. This property is interesting because it contradicts the very motivation underlying the use of non-equilibrium models. Therefore, we argue that the present version of the GSOM-like models could have minimal implications for CAV traffic flows, as investigating the instability of CAV traffic flow is a major aspect (Ciuffo et al., 2021; Gunter et al., 2021; Sun et al., 2020a).

Meanwhile, a different approach was adopted by Jiang et al. (2002) for eliminating "faster-than-traffic" waves, where another model was developed from the continuum approximation of the full velocity difference car-following model (Jiang et al., 2001). The acceleration mechanism in this model is

$$\partial V / \partial t + V\partial V / \partial x = (V_{\mathrm{e}}(\rho) - V)/\tau + C_0\partial V / \partial x \tag{3.6}$$

where $C_0$ is a constant and $\lambda_1 = V - C_0$ and $\lambda_2 = V$ are the characteristic speeds. This model is proposed as the continuum equivalent of the full velocity difference car-following model (FVDM), and in recent years, several empirical and theoretical studies have shown the capability of the FVDM for modeling AVs (Gunter et al., 2021). Therefore, at first sight, the model in Eq. (3.6) seems to be a good candidate for continuum modeling of the automated vehicles (AVs)' traffic flow.

However, this model is derived based on oversimplifying derivation processes, and the constant propagation speed $C_0$ results in several theoretical and behavioral issues. For instance, for sufficiently large relaxation times, the acceleration equation gets decoupled from the continuity equation, and when applied to real-world traffic, the model can violate the anisotropy property near the ramp areas (Mohammadian et al., 2021c). Such issues have significant implications for the CAV traffic flow. From the perspective of this model, vehicles always react to the speed differences ahead regardless of the available spacing. This aspect conflicts with how CAVs respond based on collision avoidance systems and time to collision, taking into account both speed and spacing (Van Arem et al., 2006; Xiao et al., 2017). Therefore, we argue that the present version of the model above has little implication for CAVs.

### 3.2. *Non-equilibrium models for CAV traffic flows*

This section first reviews some of the non-equilibrium models, originally proposed for human-driven traffic, that have significant implications for CAV-related issues, and it then presents a focused review of the existing non-equilibrium models developed for CAVs.

#### 3.2.1. *Models proposed for human-driven traffic that have considerable implications for CAVs*

Colombo (2002) argued that the ARZ model, presented in Section 3.1.2 does not include an explicit maximum density, and can lead to non-physical traffic states (such as negative speeds or densities above the maximum) in highly congested states. Colombo (2002) proposed another model that is known as the phase-transition model since it considers distinct hyperbolic equations for the free flow and congestion phases, respectively. For the free-flow phase, traffic is in equilibrium, and the LWR model is selected, whereas for the congestion phase, a phenomenological non-equilibrium model is proposed. The mathematical formulation of the phase transition model can be expressed as

$$\begin{cases} \partial\rho/\partial t + \partial(\rho V_{\mathrm{f}}(\rho))/\partial x = 0, \quad \text{in free flow}(\Omega_{\mathrm{f}}) \\ \begin{cases} \partial\rho/\partial t + \partial(\rho V_{\mathrm{c}}(\rho, q))/\partial x = 0, \\ \partial q/\partial t + \partial((q - q_*)V_{\mathrm{c}}(\rho, q))/\partial x = 0, \end{cases} \quad \text{in congestion } (\Omega_{\mathrm{c}}) \end{cases} \tag{3.7}$$

where $q$ is a flow-type variable, which is clarified in more detail below, $q_*$ is a constant parameter, and $\Omega_{\mathrm{f}}$ and $\Omega_{\mathrm{c}}$ are the solution domains in the free flow and congestion regions, respectively, given as

$$\begin{cases} \Omega_{\mathrm{f}} = \{(\rho, q) \in [0, \rho_{\max}] \times [0, +\infty], V_{\mathrm{f}}(\rho) \geq V_{\mathrm{f}-}, q = \rho V\} \\ \Omega_{\mathrm{c}} = \left\{ \begin{matrix} (\rho, q) \in [0, \rho_{\max}] \times [0, +\infty], \quad V_{\mathrm{c}}(\rho, q) \leq V_{\mathrm{c}+} \\ \dfrac{(q^- - q^*)}{\rho_{\max}} < \dfrac{(q - q^*)}{\rho_{\max}} < \dfrac{(q^+ - q^*)}{\rho_{\max}} \end{matrix} \right\} \end{cases} \tag{3.8}$$

where $V_{\mathrm{f}-} = V_{\mathrm{f}}(\rho_{\mathrm{cr}})$ is the minimum speed in free flow with $\rho_{\mathrm{cr}}$ being the density corresponding to the interface between the free flow and congestion regions, and $V_{\mathrm{c}+}$ is the maximum speed in congestion, and $q^-$ and $q^+$ are the minimum and the maximum values for $q$. Traffic speed in the free flow and congestion regions is given as

$$V_{\mathrm{f}}(\rho) = \left(1 - \frac{\rho}{\rho_{\mathrm{jam}}}\right)V_{\max} \text{ and } V_{\mathrm{c}}(\rho, q) = q\left(\frac{1}{\rho} - \frac{1}{\rho_{\mathrm{jam}}}\right) \tag{3.9}$$

where $V_{\max}$ is the maximum speed, $\rho_{\mathrm{jam}}$ is the maximum density under the jam condition. Characteristic speeds of the model in the congestion region are obtained as

$$\lambda_1 = (2/\rho_* - 1/\rho)(q_* - q) - q_*/\rho_* \text{and} \lambda_2 = V > \lambda_1 \tag{3.10}$$

where $\rho_*$ is another calibration parameter subject to the condition $(q/q_* + \rho/\rho_* > 1)$. The model by Colombo (2002) has some promising potential with respect to the empirical observations as well as the analytical properties of traffic flow. To elaborate, we shall interpret the quantity $q$ from a physical perspective. From micro-to-macro transforms, the average spatial gap (inter-vehicular spacing) between vehicles can be approximated as

$$S = \frac{1}{\rho} - \frac{1}{\rho_{\text{jam}}} \tag{3.11}$$

Meanwhile, by comparing Eq. (3.11) and traffic speed in the congestion phase of the phase-transition model (Eq. (3.7)), one can obtain the average spatial gap as $S = V_c/q$. This leads to obtaining the average time gap $T_g$ as $T_g = S/V_c = 1/q$, and consequently, the flow-type quantity $q$ in the phase-transition model can be interpreted physically as traffic flow corresponding to the time gap. Considering this physical interpretation, we posit that the quantity $q_*$ in Eq. (3.7) must be interpreted as the inverse of the equilibrium time gap, rather than as "characteristic of the road under consideration", the interpretation given originally by Colombo (2002). Our proposed interpretation is consistent since with $q = q_*$, the second dynamic equation in the congestion phase is eliminated, and one obtains the steady-state condition of the model in the congestion region. Under the steady-state condition, traffic flow is obtained as $Q = q_*(1 - \rho/\rho_{\text{jam}})$, which is equivalent to the triangular fundamental diagram in the congestion region. Therefore, it could be argued that under the steady-state condition, the phase-transition model by Colombo (2002) reduces to the cell-transmission model.

Blandin et al. (2011) proposed a general phase-transition model in which the triangular fundamental diagram is adopted for the free-flow phase, and for the congestion phase, the solution domain is defined as an area surrounding the triangular fundamental diagram as well. Blandin et al. (2013) investigated the performance of the phase transition model against the NGSIM data and suggested that the model can explain many traffic phenomena such as hysteresis and traffic instabilities. However, since the spatiotemporal range of the NGSIM data was quite limited and dominated by the congestion phase, the finding in Blandin et al. (2013) mainly demonstrates the model's performance regarding convective propagation of congested states. The question remains about the performance of the phase-transition models regarding aspects such as congestion emergence, queue formation, and dissolution of jam fronts.

The model by Colombo (2002) is perhaps the first non-equilibrium model that explicitly restricts the solution domain to specific ranges ($\Omega_f$ and $\Omega_c$) in the flow-density diagram. This property makes the model potentially more realistic, but there are two potential issues with the way the solution domains are introduced. First, they are discontinuous at the interface between the free flow and congestion, meaning that, the model systematically dismisses the possibility of continuous steady-state traffic from the free flow to congestion, observed by Cassidy (1998). Second, the admissible range is imposed externally and highly depends on the choice of the parameters $Q^-$ and $Q^+$. As a result, the intermediate states arising in the Riemann problems may not necessarily fall within the solution's admissible range, and imposing restrictions in the numerical solutions might be necessary, which can, in turn, result in unrealistic behavioral implications. For instance, numerical examples presented in Chalons and Goatin (2008) showed that multiple intermediate states and jumps arise for a typical deceleration corresponding to a transition from free flow to congestion.

Another line of research has utilized the gas-kinetic theory and proposed incorporating the non-local drivers' anticipation of traffic condition ahead. Treiber et al. (1999) proposed the first GKT-based non-equilibrium model in this category, where a car-following relationship also determines driver behavior at the microscopic level. Central to the underlying car-following relationship is that the traffic condition ahead ("interaction zone") also affects the steady-state speed which is derived implicitly. The acceleration equation in this model is defined as

$$\partial V/\partial t + V\partial V/\partial x = \left(V_e^*(\rho, V, \rho_a, V_a) - V\right)/\tau - 1/\rho\partial\left(\sigma_V^2(\rho, V)\right)/\partial x \tag{3.12}$$

where $\sigma^2(\rho, V) = A(\rho)V^2$ with $A(\rho) = A_0 + \Delta A[1 + \tanh((\rho - \rho_{\text{cr}})/\Delta\rho)]$ being the variance perfector in which $A_0$ and $\Delta A$ are constants, $\rho_{\text{cr}}$ is the critical density, and $\Delta\rho$ is the width of transition from free flow to congestion. $V_e^*$ is the "target speed" which is found by

$$V_e^*(\rho, V, \rho_a, V_a) = V_{\max}\left(1 - \frac{A(\rho)}{A(\rho_{\text{jam}})}\left(\frac{\rho_a T}{1 - \rho_a/\rho_{\text{jam}}}\right)^2 B(\delta V)\right) \tag{3.13}$$

where $\rho_a$ and $V_a$ are the non-local state variables at the interaction zone $(X_a = X + \gamma(1/\rho_{\text{jam}} + VT))$ in which $\gamma$ is the anticipation factor, and $B(\delta_V) = 2[1/\sqrt{2\pi}\exp(-\delta_V{}^2/2) + (1 + \delta_V{}^2)$ is the Boltzmann factor with $\delta_V = V - V_a/\sigma_V$.

Like diffusion terms, the non-local interaction in this model introduces smoothing effects on sharp shockwaves and subjects them to a finite physical width. However, the effects of non-localities are propagated only forward which prevents the model from exhibiting the wrong-way-travel effect (Helbing and Johansson, 2009). The steady-state equilibrium condition in this model is implicitly given based on the definition $\partial\rho/\partial x = \partial V/\partial x = 0$ which is obtained with $\rho = \rho_a$ and $V = V_a$. The non-local nature of speed adaptation in this model has potentials to incorporate the properties of CAVs, and Section 3.2.2 reviews the GKT-based CAV-related model.

While the implicit steady-state condition in this model is consistent with the underlying behavioral premises, the non-local GKT model is prone to two behavioral issues.

First, the model is inconsistent with the LWR model at the zero-relaxation time limit. To elaborate, we multiply the acceleration equation given in Eq. (3.13) by relaxation time ($\tau$) and take the limit $\lim \tau \rightarrow 0$, which results in $V = V_e^*(\rho, V, \rho_a, V_a)$, which does not necessarily satisfy the steady-state equilibrium condition unless $\rho = \rho_a$ and $V = V_a$. Furthermore, the function $V_e^*$ can take negative values because some parts of the non-local speed adaptation are placed in this term. The non-local GKT model can, therefore, yield non-physical results (i.e., $V = V_e^*(\rho, V, \rho_a, V_a) < 0$) at the zero-relaxation time limit.

Meanwhile, GKT-based models with two dynamic equations are prone to a conceptual inconsistency with the stead-state traffic. To elaborate, let $\rho_{\text{st}}$ and $V_{\text{st}}$ be the density and speed in the steady-state traffic, which, by definition, is any traffic condition subject to $\partial\rho_{\text{st}}/\partial t = \partial\rho_{\text{st}}/\partial x = V_{\text{st}}/\partial t = \partial V_{\text{st}}/\partial x = 0$ for a sizeable window in time-space diagram (Zhang, 1999). From the perspective of the GKT model, the steady-state traffic condition is restricted to those traffic conditions, where there exists speed variance among the underlying vehicles given as $\sigma_{V_{\text{st}}}^2 = A(\rho_{\text{st}})V_{\text{st}}^2$ subject to $\partial\sigma_{V_{\text{st}}}^2/\partial t = \partial\sigma_{V_{\text{st}}}^2/\partial x = 0$. This restriction brings two downsides. First, such a steady-state traffic condition may not frequently be observable in real-world traffic because in practice, speed variance leads to the spatial dispersion of vehicles over time violating both $\partial V_{\text{st}}/\partial t = 0$ and $\partial V_{\text{st}}/\partial x = 0$ for a moving window travelling with the average speed. Second, the steady-state condition of the GKT model is inconsistent with the car-following steady-state traffic where vehicle's speed should be invariant across different vehicles.

Nevertheless, it has been shown that the non-local GKT model can replicate many empirical observations of real-world traffic, such as the widely scattered traffic states and traffic instabilities in the flow-density diagram and scattering (Helbing, 2001; Mohammadian et al., 2021c; Treiber and Helbing, 1999). Ngoduy (2012b) extended the model by Treiber et al. (1999) to incorporate "aggressive" and "timid" driving behaviors, by assuming that aggressive drivers underestimate the downstream traffic density whereas timid drivers would overestimate the downstream traffic density. Through numerical experiments, the authors concluded that traffic flow becomes more stabilized with aggressive driving, whereas the width and amplitude of traffic instabilities increase with timid driving. The mechanisms proposed for human factors are however systematic, and do not take into account the psychological process underlying drivers' reactions to the perception of traffic. Moreover, it does not necessarily involve aggressiveness or timidness, and could be more relevant to issues such as distractions and perception of risk (Fuller, 2011; Fuller et al., 2008).

Hoogendoorn and Bovy (2000) derived another GKT-based model for multi-class traffic, where the underlying gas-kinetic model (e.g., the model by Hoogendoorn and Bovy, 1998) is a multi-class extension of the model by Paveri-Fontana (1975). Hoogendoorn and Bovy (2000) showed that the model is generic in the sense that with a proper choice of parameters, the model can be identical to some previous models (e.g., Helbing, 1996b; Payne, 1971; Phillips, 1979).

Meanwhile, Tampère et al. (2003) followed a distinct approach, focusing on incorporating human factors in the continuum framework. They derived a GKT-based model, known as the "human-kinetic" model, in which driver alertness was introduced as an additional parameter. Driver alertness in this model explicitly affects the car-following behaviors of vehicles, i.e., that vehicles tend to maintain the desired spacing and to synchronize their speed with that of their leading vehicles.

#### 3.2.2. Models proposed for CAVs and mixed CAVs

Tampère et al. (2009) proposed a multi-class variation of the human-kinetic model by Tampère (2005) for mixed traffic flows involving conventional and connected vehicles. The model assumes that conventional and connected vehicles have the same car-following properties, but the level of driver alertness is higher for connected vehicles, due to receiving information about the traffic ahead in advance. The authors suggested that traffic flow stability increased with an increase in the penetration rate of connected vehicles.

Ngoduy et al. (2009) presented another GKT-based continuum model for mixed traffic flows involving conventional and connected vehicles, where the information can travel between the connected vehicles. The model assumes that when receiving the information, connected vehicles can start decelerating smoothly before perceiving the congestion in their visible range. The propagation of information for vehicles follows the form:

$$\partial\varphi_u/\partial t + \xi\partial\varphi_u/\partial x = -\gamma\varphi_u \tag{3.14}$$

where $\varphi_u$ is a smooth normalized function in the range of [0, 1] and defines the probability of vehicles upstream receiving information ($\varphi_u = 0$ for a conventional vehicle), and $\xi$ is the propagation velocity of the message, which reflects the transmission range. The term $\gamma\varphi_u$ is introduced to nullify the effects of information beyond a certain spatial range, with $\gamma$ acting as the damping rate. The function is then incorporated into the acceleration equation at the car-following level. Ngoduy et al. (2009) derived a multi-class non-equilibrium model from the gas-kinetic model by Hoogendoorn and Bovy (2000). With the same continuity equation as in Eq. (2.6), the acceleration equation is

$$\partial V_u / \partial t + V_u\partial V_u / \partial x = (V_{e,u} - V_u) / \tau_u - 1 / \rho_u\partial P_u / \partial x \tag{3.15}$$

where $P_u = \rho_u\Theta_u$ is the so-called pressure function defined based on the velocity variance and $V_{e,u}$ is the equilibrium speed of class $u$, defined as

$$V_{e,u} = V_{\max,u} - (1 - \eta_u)\tau_u \sum_{1}^{U} \Pi_{us} \tag{3.16}$$

with $\eta_u$ being the probability of lane-changing maneuvers for class $u$ and $\Pi_{us}$ is the braking rate between vehicle class $u$ and leading vehicle of any class. Ngoduy et al. (2009) used the model to study the mixed traffic flow for different traffic compositions, and the results suggest that both the capacity and stability of traffic flow improve as the penetration rate of connected vehicles increases. Ngoduy (2012a) presented another continuum model for mixed-traffic flows of conventional and automated vehicles from the gas-kinetic model by Hoogendoorn and Bovy (2000), in which the acceleration mechanisms of ACC-equipped vehicles obey the car-following relationship provided by Davis (2004). The acceleration equation for this model is

$$\partial V_u / \partial t + V_u\partial V_u / \partial x = (V_{e,u} - V_u) / \tau_u - 1 / \rho_u\partial P_u / \partial x + (\gamma_u / \tau^*)\partial V / \partial x \tag{3.17}$$

where $\gamma_u = 1$ for ACC vehicles, and $\gamma_u = 0$ for conventional vehicles. Similar to previous works, the numerical tests for this model also suggest that the stability and the capacity of traffic flow improve with increases in the penetration rate of ACC-equipped vehicles.

Ngoduy (2013a) extended the model by Treiber et al. (1999) for studying traffic flows of vehicles equipped with cooperative adaptive cruise control (CACC) driving strategies, where the CACC element is derived from a multi-anticipative car-following relationship. The acceleration equation in this model is presented based on traffic flow as

$$\partial V / \partial t + V\partial V / \partial x = \left(V_e^*(\rho, V, \rho_a, V_a) - V\right) / \tau - 1 / \rho\partial\left(AV^2\right) / \partial x + \nu_{ACC} \tag{3.18}$$

where $\nu_{ACC} = 1/\rho\partial V/\partial x \sum_{n=1}^{N} n/\tau_n^*$ incorporates the acceleration of Cooperative Adaptive Cruise Control (CACC) vehicles, with $N$ being the number of preceding vehicles that a follower can interact with, and $N = 1$ accounting for ACC traffic flow without cooperation. $\tau_n^*$ is the state-dependent relaxation time of vehicles to adapt their speed to their $n^{th}$ leader and is defined as

$$1 / \tau_n^* = 1 / (2\tau_m^0)[1 + \tanh((\rho - \rho_{cr}) / \Delta\rho)] \tag{3.19}$$

where $\tau_m^0$ is the maximum relaxation time of vehicles under the jam condition. It is assumed that the condition $\tau_1^* < \tau_2^* < \ldots < \tau_N^*$ is satisfied within a platoon. Numerical simulations and stability analysis suggested that both CACC and ACC traffic have greater flow capacity and stability range than conventional traffic, with CACC traffic performing the best. Another continuum model for automated vehicle platoon-based driving was presented by Ngoduy (2013b). The model is based on the GKT-based model by Hoogendoorn and Bovy (2000), where the car-following model by Helbing and Tilch (1998) was extended for a platoon of vehicles. It was established in the aforementioned study that platoon-based driving significantly increases the stability of traffic.

Delis et al. (2015) proposed a generic variation of the GKT-based model for comparing the traffic flow of ACC and CACC vehicles with that of conventional traffic. The acceleration equation in this model is

$$\partial V / \partial t + V\partial V / \partial x = \left(V_e^*(\rho, V, \rho_a, V_a) - V\right) / \tau(1 - \beta F(\rho)) + \left(\alpha\upsilon_{ACC} - 1 / \rho\partial\left(AV^2\right) / \partial x\right) \tag{3.20}$$

where $\alpha = \beta = 1$ defines CACC driving, $\alpha = 1$ and $\beta = 0$ define ACC driving, and $\alpha = \beta = 0$ describes conventional human-driven traffic as in Treiber et al. (1999), where $F(\rho) = 1 + \tanh((\rho - \rho_{cr}) / \Delta\rho)$ and the term $\upsilon_{ACC}$ is derived by

$$\nu_{ACC} = \frac{1}{\rho} 0.5 \sum_{n=1}^{N} \frac{\rho^* V_n^* - \rho V}{\tau_n^*} F(\rho) \tag{3.21}$$

where $V_n^*$ refers to the non-local speed for vehicle $n$ computed at $X_n^* = X + \gamma^*(1 / \rho_{jam} + T^* V)$, and $\rho^* = 1 / (1 / \rho_{jam} + T^* V_n^*)$ and other parameters are previously defined. While the model by Ngoduy (2013a) defines the acceleration of ACC in terms of gradients, Delis et al. (2015) adopted a non-local approach and incorporated the acceleration of ACC/CACCs into the relaxation term. Similar to the work of Ngoduy (2013a), the study concluded that both ACC and CACC vehicles improve the stability of traffic flow. Compared to the model by Ngoduy (2013a), the stabilizing effects of ACC/CACC vehicles were more significant in the model by Delis et al. (2015). The instability range of the model is analytically discussed in more detail in Porfyri et al. (2015).

Delis et al. (2018) proposed a multi-lane variation of the model, where each lane has its own dynamic equation, allowing vehicles to change lanes depending on the relative traffic condition between the lanes. Such an approach may be suitable for CAV traffic flow as it can address some of the issues of continuum modeling associated with the single-pipe treatment of traffic flow. However, the model by Delis et al. (2018) may become intricate for freeways with many lanes. In this model, the mechanism for lane-changing maneuvers in this model is phenomenological and based entirely on the relative traffic condition between the adjacent lanes, and does not account for the human factors and decision-making processes involved in the adequately reflect the human factors of lane-changing maneuvers and their effects on traffic flow properties (Zheng, 2014).

Zheng et al. (2015) derived a non-equilibrium model accounting for a situation in which drivers receive information about traffic conditions behind, with potential applicability to the connected environment. The resulting acceleration equation in this model is derived from the car-following found by Helly (1959).

$$\partial V / \partial t + (V - C_0)\partial V / \partial x = r_f \alpha_f \left(1 / \rho - 1 / \rho_e^f(V)\right) - r_b \alpha_b \left(1 / \rho - 1 / \rho_e^b(V)\right) \tag{3.22}$$

where $C_0$ is a constant sonic velocity, indices $f$ and $b$ refer to traffic condition ahead and behind respectively, $\rho_e(V)$ defines the equilibrium density, $r$ is a weight factor defined for balancing the attention to traffic condition ahead and behind, and $\alpha(1/\mathrm{s}^2)$ is the sensitivity coefficient which is comparable to relaxation time. In this model, the effects of the available information on the acceleration are placed in source terms based on the local traffic condition. However, the available information on both sides is expected to affect the anticipatory acceleration in real-world traffic. Ngoduy and Jia (2017) incorporated the multi-anticipative driving in the car-following model by Helly (1959) and derived a generic variation of the model by Zheng et al. (2015), where the effects of multi-anticipative driving are placed into a gradient term. The authors suggested that multi-anticipative driving improved the operational properties and the stability of traffic flow.

Overall, the numerical simulation outcomes from the existing models reviewed in this section have consistently suggested that CAVs deliver positive effects for traffic flow properties such as increasing traffic capacity, improving traffic flow stability, and reducing congestion emergence and propagation. However, these outcomes cannot be generalized to the real-world traffic flows of CAVs for several reasons:

(1) Existing CAV models are not flexible enough to accommodate a wide range of potential impacts of different CAV penetration rates on traffic flow dynamics. Almost all the existing studies on continuum modeling of CAVs suggest that with an increase in their penetration rate, CAVs can systematically bring about positive outcomes (e.g., increased capacity) for traffic flow, and their stabilizing effect can emerge at small penetration rates (less than 30%) regardless of how CAVs are spatially distributed in the mixed traffic (e.g., Ngoduy et al., 2009). On the other hand, some recent microscopic simulations have suggested that there could be a significantly large threshold (e.g., 70%) for the penetration rates of CAVs, below which automated driving may not deliver significant improvements to traffic flow (Calvert et al., 2017).

(2) Existing models dismiss the possibility of negative effects of AVs on traffic flow properties, which have been observed in some field experiments. For instance, several recent empirical observations of ACC-equipped vehicles have suggested that such vehicles are string unstable (Ciuffo et al., 2021; Gunter et al., 2021). Recent empirical analysis of Waymo vehicles' trajectories have revealed that CAVs drive more conservatively for safety considerations, and such driving strategies decrease traffic efficiency (Hu et al., 2023). However, in many existing continuum models (Delis et al., 2015, 2018; Ngoduy, 2013a, 2013b; Ngoduy et al., 2009; Ngoduy and Jia, 2017; Ngoduy and Wilson, 2014), CAVs systematically stabilize traffic flow, which could be, in part, caused by the discrepancy between systematic speed adaptation in continuum models and CAVs' hierarchical control at the car-following level. Section 4.1 delves into this issue in more detail.

(3) Finally, existing models are based on the notion that drivers fully trust in automated vehicles (AV) and avoid interventions with AVs' driving, which contradicts the empirical findings, suggesting that drivers are likely to intervene with AVs depending on traffic conditions and the type of car-following situations. For instance, through field experiments with AVs, Viti et al. (2008) found that drivers tend to deactivate ACC driving in congested traffic conditions, which are important contexts for the efficacy of the ACC system. Varotto et al. (2015) conducted a driving-simulator study and found that such interferences with AVs negatively affected traffic flow properties (e.g., resulting in increased time gap), and thereby could undermine the expected positive effects of CAVs. Recent field experiments with AVs and driver behavior modeling have suggested that drivers' perceived risk and task difficulty level are major determining factors with respect to driver's intervention with automated vehicles and transitions from manual to automated driving (Varotto et al., 2017, 2018).

## 4. Discussion: continuum modeling of real-world CAV traffic flows

This section discusses some challenges and opportunities inherent to commercial CAVs and revisits their implications for continuum modeling, while highlighting some major problems and proposing potential approaches for addressing them. Our focus is on continuum modeling of *real-world* CAV traffic flows, once commercial CAVs have been introduced in traffic streams.

Our discussions are organized into two parts. First, we focus on the properties, challenges, and opportunities inherent to commercial CAVs, and revisit their implications on continuum modeling. Next, we discuss and highlight some essential challenges and problems pertaining to continuum approximation of real-world (mixed) CAV traffic and propose some potential solutions.

### *4.1. Some unique properties, challenges, and opportunities of CAVs and their implications for continuum modeling of real-world traffic flows*

CAVs have three powerful resources that if utilized cooperatively and with their maximum potentials, can lead to unprecedented paradigm shifts in traffic flow modeling and operations (Shladover, 2018). These resources are:

(1) High-fidelity monitoring sensors such as cameras equipped with image-processing power, radar, lidar, etc. that collect high-resolution information in their detection range. Some examples include spatial gap, speed difference between their followers and leaders, the number of leaders ahead in their current and adjacent lanes, etc.

(2) Communication systems that enable them to exchange a wide range of information with other CAVs and infrastructures, including information about traffic conditions beyond their sensors' detection ranges, warnings about congestion and accidents, specific driving rules, desirable driving objectives (e.g., eco-driving, cooperative driving).

(3) Programmable AI-powered computing systems that can quickly process the information provided by (1) and (2), define, adapt, and execute the driving strategy accordingly in real time.

Motivated by maximizing the use of such resources, numerous studies have developed optimal driving strategies for CAVs with respect to certain objectives, such as improving traffic capacity (Guanetti et al., 2018; Rios-Torres and Malikopoulos, 2018), ensuring string and traffic flow stability (Yao et al., 2020; Zheng et al., 2020), alleviating fuel-consumption (Taiebat et al., 2018; Yan and Wu, 2021), and minimizing crash risk (e.g., Wang et al., 2020). The existing continuum models developed for CAVs reviewed in the previous sections also fall into this category.

Meanwhile, some studies have focused on the use of individual or multiple CAVs as moving controllers to improve traffic flow properties, such as: (a) congestion alleviation, e.g., by proposing shockwave mitigation strategies (Talebpour et al., 2018), (b) dissipation of stop-and-go traffic and stabilizing traffic flow (Sun and Tan, 2020), and (c) enhancing traffic safety by proposing driving strategies that would alleviate crash risk (Li et al., 2017; Morando et al., 2018), where some recent studies have integrated both traffic efficiency and safety considerations in their controllers (e.g., Han et al., 2021b).

In short, such studies focus on how CAVs should operate ideally in order to achieve certain objectives. Obviously, this is an important research direction. However, this type of research may not necessarily provide sufficient understanding about the actual properties of CAVs and mixed CAVs' traffic flows at the macroscopic level once they are actually in place. Another point is that the objectives considered from safety, environmental and traffic flow perspectives are not necessarily aligned and in some cases are in sharp conflict with each other (Li, 2022; Mohammadian et al., 2021a).

Meanwhile, since CAVs have not yet been widely introduced in traffic, there are many uncertainties in their expected potential benefits, as highlighted by some conflicting outcomes regarding the impacts of CAVs on traffic flow dynamics drawn from the limited existing field tests and empirical observations of CAV traffic flows (Ciuffo et al., 2021; Stern et al., 2018).

For instance, in some cases, the findings from field testing suggest that CAVs can stabilize traffic flow, and dissipate the instabilities caused by human-driven vehicles (Stern et al., 2018), whereas in other cases, it has been shown that some commercial CAVs are not string stable (Ciuffo et al., 2021). In addition, studies on the field tests with Waymo data (Hu et al., 2022; Sun et al., 2020a) has revealed that Waymo vehicles have more conservative gap acceptance behaviors compared to human-driven vehicles (Hu et al., 2023), and such findings are in sharp contrast with the assumptions underlying many models developed for CAVs, where properties such as decreased time-gap and platooning of CAVs had been utilized to increase traffic flow capacity (Schoenmakers et al., 2021; Sun et al., 2020a; Tilg et al., 2018; Wu et al., 2022).

The remaining parts of this section discuss some properties, opportunities, and challenges of commercial CAVs, using established knowledge as well as empirical findings in the literature. It then revisits their implications on continuum modeling of real-world (mixed) CAV traffic flows.

#### *4.1.1. CAVs' driving strategies, drivers' interactions with CAVs, and the pronounced role of human factors*

Due to their AI-powered detection sensors and powerful computing resources, CAVs will use complex car-following driving strategies that might differ, in many ways, from the analytic automated driving rules in the majority of existing car-following models (Guanetti et al., 2018; Guvenc et al., 2017).

At the same time, commercial CAVs' exact driving strategies will remain a puzzle in the foreseeable future due to factors such as market competition and intellectual property concerns, which make manufacturing companies that produce CAVs less likely to release the exact algorithms and procedures utilized to program their CAVs' responses and maneuvers in different situations. It is also reasonable to assume that different CAV companies will design the driving strategies differently considering numerous factors, ranging from vehicle's mechanical properties, state laws and design guidelines for CAVs, etc., to real-time situational decision-making criteria (e.g., safety considerations (Zhao et al., 2021), drivers' comfort, drivers' preferred driving setting). In other words, even if CAV traffic is operated at the highest automation level and human factors are removed, there would still be heterogeneities among CAVs regarding driving signature.

Such uncertainties would pose additional challenges when it comes to investigating the impacts of CAVs' underlying control laws on traffic flows' collective properties, such as wave propagation, stability, etc. For instance, consider a scenario where CAVs' longitudinal control laws to determine car-following rules are based primarily on safety considerations. If minimizing rear-end crash risk is the key criterion for identifying the optimal response to the situation ahead, one likely scenario is that CAVs would use time-to-collision in each time step as a measure of crash risk, and update their position for the next time step such that time-exposed time-to-collision (TET) is minimized during the entire deceleration maneuver using projected values of time gap and speed values for each possible move (Dai et al., 2023; Yu et al., 2021). The implications of such safety-oriented responses may not necessarily be desirable, from a traffic flow dynamics perspective, in terms of aspects such as the waves initiated by the aforementioned vehicle and their effects on the upstream vehicles.

CAVs are also likely to bring about new dimensions to human psychological factors that underly drivers' perceptions and reactions to car-following situations, which will ultimately affect the collective properties of traffic flow. A relevant example is the situation in which the system may request that drivers take over the vehicle due to uncertain situations ahead (poor lane-marking), system errors, etc. (Zhang et al., 2019). Many driving simulator studies and empirical studies of real-world observations of CAVs have found that human psychological factors, including the mental workload preceding the take-over maneuvers and drivers' perception of the situation (regarding risk and task difficulty) would affect both drivers' take-over time and the quality of their take-over maneuvers (Bueno et al., 2016; Ko and Ji, 2018; Naujoks et al., 2018; Varotto et al., 2017, 2018; Zhang et al., 2019). While a detailed understanding of drivers' take-over performance is not relevant to continuum modeling, the time taken by drivers to take-over the vehicle (i.e., take-over time) as well as the quality of their take-over performance can have significant consequences on the traffic flow.

From a traffic flow theory perspective, drivers' take-over time can be considered to be the response time to a stimulus (i.e., take-over request), in which case, a stimulus-response framework can be utilized within the car-following theories (Saifuzzaman and Zheng, 2014) to investigate the consequences on the car-following behaviors of the subject vehicle and its followers. Theoretical and empirical studies of car-following behaviors have shown that drivers' response time is a major determining factor for many traffic phenomena, such as emergence and propagation of kinematic waves and traffic instabilities (Chen et al., 2012; Kesting and Treiber, 2008; Zhang and Kim, 2005; Zheng et al., 2011a).

These points may also apply when CAVs are used in the human-driven mode but provide the drivers with information about the surrounding traffic condition outside the drivers' perception zone. In such a condition, the potential impacts of the connectivity on traffic flow will, again, depend primarily on how drivers process and react to the information. Recent driving-simulator studies on driving in connected environments have found that drivers' response time may be increased in certain car-following situations (Sharma et al., 2019c), and that such behavioral adaptations can be explained with psychological theories of driver behavior (Sharma et al., 2019b, 2020). Furthermore, recent empirical and theoretical studies have underscored profound linkages between drivers' response time and human psychological factors such as task-difficulty, driver compliance, and risk perception (Saifuzzaman et al., 2015).

While a detailed understanding and modeling of driver behavior and human psychological factors are not the focus of continuum modeling, such human-factor related issues at the individual level can initiate a chain of car-following interactions that could become quite relevant for macroscopic traffic flow dynamics, such as new or unexpected patterns of wave propagations and instabilities. For instance, a single take-over maneuver in a line of automated vehicles may, depending on drivers' take-over time and quality, result in a fast backward propagation of shockwaves.

One way to investigate these issues would be to develop hybrid continuum models that incorporate both automated and manual driving for the class of automated vehicles. Such hybrid models must distinctively incorporate robust and realistic speed adaptation mechanisms for both automated and manual driving phases. The most critical component of such a hybrid modeling approach would be designing transition mechanisms from automated to manual driving and vice versa. For example, robust mechanisms must consider major human psychological factors involved in such transitions (e.g., drivers' boredom, risk perception.) as well as their impacts on driver behavior (e.g., take-over time).

#### 4.1.2. CAVs' fundamental diagrams and the implications for continuum modeling of (mixed) CAV traffic flows

In terms of modeling and empirical investigation of (mixed) CAV traffic flows, two important research questions have probably attracted the most attention in the literature: (1) How would CAVs affect macroscopic aspects such as traffic capacity, and patterns of congestion evolution, propagation, and dissipation? And (2) how would CAVs affect the stability of traffic flow? In this section, we primarily focus on the former and highlight some essential problems by considering the different components relevant to fundamental diagrams of CAVs and their implications for continuum modeling.

Our discussions primarily concentrate on the properties of CAVs, especially at the car-following level. To support this, we utilize empirical observations and field tests from automated driving, specifically with commercial adaptive cruise control vehicles (e.g., Ciuffo et al., 2021; Gunter et al., 2021; Li et al., 2021, 2022b; Shi and Li, 2021). We also incorporate established knowledge about CAVs and have organized our insights into distinct remarks.

> Remark (1): CAVs offer valuable trajectory data, reinforcing the existence of the fundamental diagram. However, its shape for CAV traffic might not always be triangular and can be influenced by the steady-state policies of the CAVs.

The controversies around the concept of the fundamental diagram, as deterministic bi-variate relations among macroscopic state variables, have been ongoing for decades (Banks, 1991; Cassidy, 1998; Daganzo et al., 1999; Kerner, 2016; Li and Zhang, 2011; Schönhof and Helbing, 2009; Treiber et al., 2010; Zhang, 2003). Some recent studies have shown the existence of the triangular fundamental diagram in the steady-state traffic using empirical trajectories collected from field testing and other observations of automated vehicles (Li et al., 2022b; Shi and Li, 2021). The outcomes of these studies suggest that steady-state CAV traffic flows are in much better agreement with the triangular fundamental diagram for different headway settings, where the maximum traffic capacity has been observed in the minimum headway setting. While such findings provide reliable empirical support for the concept of the fundamental diagram for CAVs, several points must be considered in terms of conclusions made from such studies, which we shall discuss in the following remarks.

Nevertheless, the fundamental diagram for CAV traffic may vary based on the steady-state policies used for CAVs, rather than being triangular. From a microscopic perspective, the triangular fundamental diagram is equivalent to constant time gap (or constant spacing) policy (Daganzo, 2006a). Therefore, the empirical support for the existence of the triangular fundamental diagram for the steady-state CAV traffic flows provided by the recent studies (e.g., Shi and Li, 2021) imply that the CAVs used in those studies most likely used constant time gap policy for gap control in the steady state. The triangular fundamental diagram might not fit real-world CAV traffic, especially if CAVs use variable spacing policies, leading to a non-linear steady-state diagram.

> Remark (2): Uncertainties surround CAVs' impact on flow-density patterns, with human factors remaining influential, and the fundamental diagram derived from CAV trajectories being inadequate for investigating the kinematic waves in mixed CAV traffic.

On this point, one could derive different conclusions, depending on how the problem is approached. From a theoretical perspective, a recent study has suggested that with an increase in their market penetration rate, CAVs will decrease the observed scattering phenomenon in the flow–density diagram (Zhou and Zhu, 2020). From a controlled empirical perspective, the outcomes of the recent studies imply that the range of scattering would be considerably limited if one develops case-specific fundamental diagrams from the macroscopic observations that match different headway settings in recent studies.

However, when it comes to the real-world CAV traffic, the scattering phenomenon might still be observed even for purely automated traffic in the steady state, due to drivers' preferred driving settings in automated driving mode and/or different control laws of CAVs from different manufacturers. For instance, in a line of automated vehicles supervised by human drivers, different drivers may use different headway settings according to their perception of risk and other human factors. From a macroscopic perspective, the variances in headway settings would show up as a range of traffic densities for a given speed. Thus, we argue that at lower levels of automation and when supervised by human drivers, CAVs' impacts on the fundamental diagram might still be at the mercy of behavioral heterogeneities between human supervisors of CAVs.

As such, the fundamental diagram, derived from the steady-state traffic is likely to be insufficient to investigate the kinematic waves in (mixed) CAV traffic flows. To elaborate, let us first review some basic aspects relevant to the concept of fundamental diagram and its linkage to the kinematic wave theory and the equilibrium continuum models:

(1) Fundamental diagram describes the relationships between macroscopic state variables (as $V_{\mathrm{e}} = V_{\mathrm{e}}(\rho)$ and $Q_{\mathrm{e}} = \rho V_{\mathrm{e}}(\rho)$), which is empirically valid only in steady-state traffic (Cassidy, 1998; Li et al., 2022b; Seo et al., 2019).
(2) For a well-defined functional form with respect to empirical observations, the slope of the fundamental diagram in the $\rho - Q_{\mathrm{e}}(\rho)$ plane is a speed-like quantity ($V^{*}(\rho) = Q_{\mathrm{e}}^{'}(\rho)$ for reference) that is restricted to a certain range $[-|C_{\mathrm{jam}}|, V_{\mathrm{cr}}]$ (Del Castillo and Benítez, 1995a).
(3) In the equilibrium continuum models, (1) and (2) are integrated with the kinematic wave theory and the propagation velocity of smooth density variations (characteristic speed $(\lambda)$) is *assumed* to be equivalent to $\lambda = V^{*}(\rho)$ .
(4) From (3), one obtains a microscopic interpretation that driver's response time to the stimulus ahead depends only on spacing and is obtained as $\tau_{\mathrm{rsp}}(\rho) = f(V^{*}(\rho))$ (Daganzo, 2006a).
(5) Empirical observations have suggested that beyond a certain congested state ($\rho_{\mathrm{cr}}^{*}$ for reference) the propagation velocity $(\lambda)$ is invariant to the density (i.e., $\lambda^{'}(\rho) = 0$ for $\rho \in [\rho_{\mathrm{cr}}^{*}, \rho_{\mathrm{jam}}]$) (Del Castillo et al., 1994; Treiber and Kesting, 2013 (chapter 8)). As such, a functional form of the fundamental diagram that gives $V^{*} = -|C_{\mathrm{jam}}|$ for $\rho \in [\rho_{\mathrm{cr}}^{*}, \rho_{\mathrm{jam}}]$ (as in the triangular-like fundamental diagrams) is desirable for equilibrium modeling frameworks (Del Castillo, 2012; Del Castillo and Benítez, 1995a, 1995b).
(6) From a microscopic perspective, if (5) holds, then driver's or vehicle's response time ($\tau_{\mathrm{rsp}}(\rho)$) would be the same as time gap ($\tau_{\mathrm{gap}}$) (which is equivalent to headway setting in CAVs) and constant over the density range $[\rho^*_{\text{cr}}, \rho_{\text{jam}}]$ (i.e., $\tau_{\text{rsp}} = \tau_{\text{gap}} = cte$) (Del Castillo et al., 1994).

(7) From (5) and (6), the equilibrium modeling frameworks suggest that any transition between traffic states in the density range $\rho \in [\rho^*_{\text{cr}}, \rho_{\text{jam}}]$ is a shockwave (for both acceleration and deceleration) maneuvers, of the same magnitude as propagation velocity as $|\lambda| = |C_{\text{jam}}| = (\tau_{\text{gap}}\rho_{\text{jam}})^{-1} = (\tau_{\text{rsp}}\rho_{\text{jam}})^{-1}$ (Zhang, 2001).

We now revisit the points above and discuss their implications for the traffic flow of CAVs.

Recall that the assumption $\lambda = \lambda(\rho) = V^*(\rho)$ in point (3) is an *ansatz* underlying equilibrium continuum models and is not a direct deduction from (1) and (2). To obtain the propagation velocity $\lambda$ in real-world, one would have to track disturbances in the density profiles in the time–space diagrams (e.g., by investigating the cross-correlation of speed time-series in different locations (Coifman and Wang, 2005; Treiber and Kesting, 2012; Zielke et al., 2008)), whereas $V^*(\rho)$ is simply the slope of fundamental diagram in the $\rho - Q$ plane, derived from the empirical observations of steady-state traffic.

By dismissing the distinction, several empirical studies of fundamental diagrams have utilized the observed slopes ($V^*(\rho)$) in order to associate them with kinematic waves (Li and Zhang, 2011; Li and Chen, 2017; Li et al., 2022b; Lu et al., 2009), where in some cases, further attempts have been made to provide empirical support for some unusual phenomena such as deceleration fans (Li and Zhang, 2011).

Likewise, recent empirical studies of fundamental diagram of commercial CAVs have also dismissed the distinction discussed above, and thus, have utilized the triangular fundamental diagram, derived from empirical CAV trajectories, to investigate the kinematic waves of CAVs. For instance, Li et al. (2022b) utilized the observed slope ($V^*(\rho)$) in the congestion branch, and concluded that the magnitude of wave speeds $|C_{\text{jam}}|$ for CAVs' shortest headway settings can be unprecedentedly large (around 100 km/hr ) and raised concerns around the safety of CAVs (Li et al., 2022b). Meanwhile, such conclusions are not in line with empirical investigations, suggesting that CAVs’ response time is generally larger and safer than average human-driven vehicles (Hu et al., 2023), at least at the current stage, meaning that the corresponding wave speeds might be even smaller in magnitude than those issued by the corresponding human-driven vehicles.

These discrepancies can be explained by several factors. First, CAVs' response time most likely depends on the nature of the car-following situations (acceleration vs. deceleration) and safety considerations such as time-to-collision (Calvert et al., 2018; Delle Monache et al., 2019; Shladover, 2018). This suggests that at least three components (i.e., the CAV's speed and spatial gap as well as the sign of speed difference between the CAV and its leader) are likely to be continuously processed in real time by CAVs in order to calculate the necessary response at each time step. As such, CAVs' corresponding wave speeds when responding to different stimuli will depend on these components and vary considerably between the acceleration and deceleration situations. Thus, the assumption that response time equals time gap is overly inadequate for modeling CAVs.

This limitation could also explain why many existing studies of CAVs developing equilibrium models have concluded that CAVs will help alleviate congestions (e.g., Jin et al., 2022; Qin et al., 2021; Wu et al., 2022). Such studies generally assumed that CAVs have shorter time gaps than human-driven vehicles, where triangular-like fundamental diagrams have been utilized to incorporate CAVs' different headway settings and their speed adaptations (see e.g., Fig. 1). In such a modeling framework, CAVs’ response time is considered the same as time gap, and with respect to point (7), under shorter time gaps, their wave speed magnitude when leaving the congestion increases, leading to a faster dissipation of congestion.

The same point applies to the findings for mixed CAV traffic within the continuum framework (e.g., Qin et al., 2021), where with an increase in their penetration rate, the average time gap in the mixed traffic is likely to increase to ensure safety and point (7) applies regarding the impacts of CAVs on congestion alleviations. However, for mixed traffic flows the situation will be much more complex because the corresponding waves could be largely dependent on how CAVs are spatially distributed between human-driven vehicles and the interactions between the two classes.

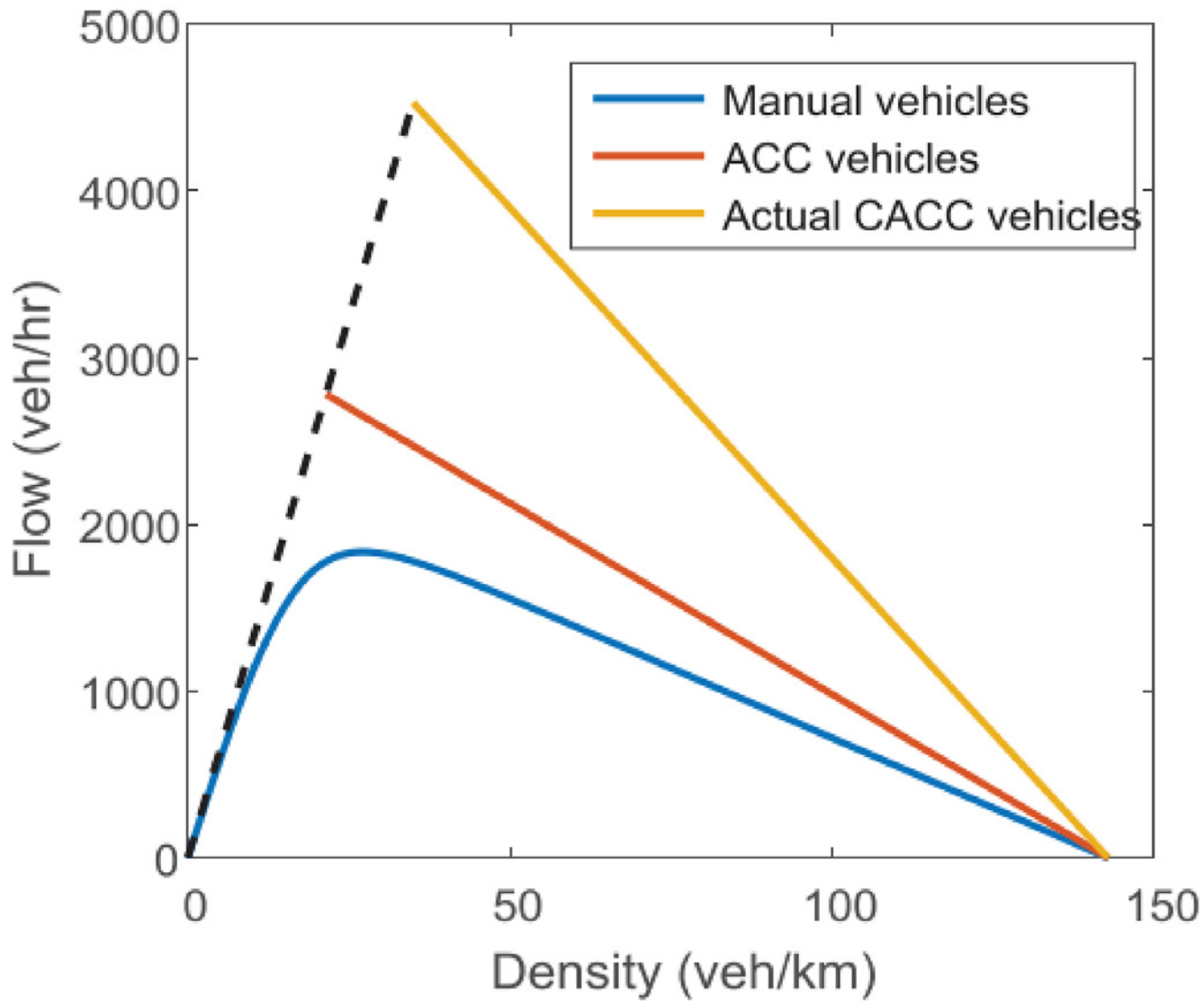

**Fig. 1.** Example of how CAVs' fundamental diagrams have been designed in several existing equilibrium continuum models. Reproduced with permission from Qin et al., © INFORMS 2021.

### 4.2. *Essential requirements, challenges and research needs for continuum modeling of (mixed) CAV traffic flows: some problems and potential solutions*

#### 4.2.1. *Need for adequately incorporating steady-state for continuum models of CAV traffic flows*

We recall from the review sections that maintaining certain car-following gaps in steady-state traffic has been a postulation underlying many car-following and continuum models of conventional and CAV traffic flow. This postulation leads to the emergence of relaxation terms in the non-equilibrium models. This section takes a closer look at the existing non-equilibrium models and discusses some considerations, requirements and challenges surrounding the relaxation terms for continuum models of CAVs.

While the majority of car-following-based and GKT-based non-equilibrium continuum models, reviewed in Section 3, employ relaxation terms, ever since the models by Aw and Rascle (2000) and Zhang (2002), there has been a distinct class of continuum models (typically those without faster-than-traffic characteristic speeds), where the relaxation term is not present in the original form of the models (e.g., Blandin et al., 2011; Garavello and Goatin, 2011; Goatin, 2006). Whether relaxation terms should be present or not in a continuum model is a dilemma of the relaxation term as there are advantages and disadvantages to both cases.

The absence of the relaxation term brings about several benefits from mathematical and computational perspectives, e.g., the models are likely to be well-posed provided that characteristic speeds are distinct over the solution domain (LeVeque, 2002). In this case, all possible Riemann problems that can arise at the discontinuities can be studied analytically, and the Godunov-type schemes can be developed to prove the existence of the unique numerical solutions Chalons and Goatin (2008); Mammar et al. (2009). This property can also be significantly advantageous for modeling traffic flow on networks because the Godunov scheme can be viewed as a supply-demand analysis (Lebacque et al., 2005), and therefore, it can be coupled with node models to determine the exact distribution of traffic flow at junctions (e.g., Garavello and Piccoli, 2006a, 2006b).

However, several conceptual problems arise in the absence of relaxation terms, which can ultimately cast doubt on whether such models should be considered non-equilibrium, and on whether they are desirable for non-equilibrium analysis of traffic flow regarding some aspects such as stability analysis.

First, in such models, any initial condition satisfying $\partial V/\partial x = 0$ is a steady-state condition, suggesting that there is no unique fundamental diagram in the models devoid of a relaxation term. To resolve this issue, Lebacque et al. (2007b) argued that the models by Aw and Rascle (2000) and Zhang (2002) can be interpreted as an extension of the LWR model with a "variable fundamental diagram". However, such an interpretation does not completely reconcile the concept of the fundamental diagram with such models. This is because another consequence of the absence of the relaxation term is that any initial condition satisfying $\partial V/\partial x \approx 0$ results in a long-lasting travelling wave.

Second, related to the first point, deviations from any initially steady-state condition may grow in amplitude and span in the time–space diagram, and never be damped. However, empirical observations and analytical studies have posited that traffic instabilities can only occur in certain intermediately congested states (Treiber and Kesting, 2011; Zheng et al., 2011a, 2011b), and that not all deviations from steady-state traffic can grow in magnitude and span in the time–space diagram (Helbing et al., 1999, 2009), i.e., under stable traffic regimes, small disturbances will diminish and converge to the steady-state equilibrium traffic. Therefore, it could be argued that models devoid of a relaxation term cannot consistently replicate the non-equilibrium traffic and the associated phenomena.

In the presence of the relaxation term, in contrast, one obtains unique (or headway-specific) fundamental diagrams, and the conceptual problems discussed above do not apply. However, an important issue is how to formulate and measure the relaxation time, whether it is constant across entire traffic regimes, and to what extent the choice of the relaxation time in such models can affect the model's capability in replicating traffic instabilities. In many continuum models with a relaxation term (reviewed in Section 3), relaxation time is treated as a constant parameter, meaning that relaxation acceleration is only dependent on the difference from the equilibrium condition and is insensitive to the traffic regime.

When it comes to CAV traffic flows, it seems that (a) utilizing the relaxation terms will be a necessity for traffic flows of CAVs, and (b) the way the relaxation term appears in the current continuum models, i.e., systematic presence in the speed-adaptation mechanisms, would be inadequate for CAVs and needs to be revisited.

Regarding (a), Section 4.1.1 discussed that there are strong empirical and analytical evidence underscoring that CAVs utilize gap-control policies in the steady state traffic, endorsing the concept of the fundamental diagram, and the relaxation phenomenon in the continuum models of freeway traffic flows. In addition, the headway control mechanism of adaptive cruise control is based on a desired headway set by drivers, meaning that automated vehicles adjust their headway even when there is no speed difference between the AVs and their leader. Relevant to this point, a key motivation behind continuum modeling of mixed CAV flows, both in the existing literature and future directions, is to understand how CAVs would affect the stability of traffic flow (Sun et al., 2018). As discussed above, this aspect can only be investigated if the models adequately incorporate the steady-state traffic, and in the absence of relaxation terms, small disturbances will not be damped, and one cannot differentiate between the stable and unstable traffic regimes (Treiber and Kesting, 2013 (chapter 15)).

Regarding (b), we discussed in Section 4.1.1 that CAVs' control laws are most likely comprised of hierarchical rules designed with respect to a wide range of considerations (e.g., the nature of the driving regime (acceleration, deceleration, safety-critical, etc.)). Therefore, gap control with respect to drivers' preferred headway is not the only criteria in the speed adaptation mechanism and may be inactive in certain situations. Perhaps a way to reconcile such relaxation terms could be to propose mechanisms to restrict the presence of the relaxation terms in certain situations.

#### 4.2.2. *Continuum model modeling of real-world (mixed) CAV flows: some fundamental challenges and possible solutions*

Section 4.14.1.2 elaborated on some complexities of the real-world mixed CAV traffic flows and discussed correspondingly the limitations of existing continuum models developed for CAV traffic. We concluded that despite their advantages, the existing models could be insufficient for real-world CAV traffic flows, from the perspectives of either understanding their macroscopic traffic flow dynamics or predictions of their operational measures for control purposes. As well, in Section 4.1.1, we discussed some fundamental challenges regarding the control laws of commercial CAVs, from the perspective of continuum modeling. In particular, we discussed the importance of investigating how the new properties brought about by CAVs' control laws and human factors could affect macroscopic traffic flow dynamics, and thus, the importance of developing new continuum models for CAV traffic. Nevertheless, the intricacies surrounding continuum approximations of CAVs' AI-powered, and discrete-like control laws can make this aim almost infeasible. This section revisits the importance of developing continuum models for CAV traffic flows and discusses some potential remedies for overcoming some existing inherent challenges.

Before delving into the technical details and challenges, it necessary to clarify the aim of model development, i.e., what would be the purpose of model application, since this is the first step in examining whether the existing challenges are relevant at all, and if so to what extent, and whether there are ways to overcome the barriers. In this regard, we distinguish between two main research directions as follows:

(a) *Developing models for explanatory purposes*, i.e., understanding the physical and analytical properties of real-world (mixed) CAV traffic flows primarily for making reproducible inferences.
(b) *Developing models for prediction purposes*, i.e., predicting the observed real-world patterns of (mixed) CAV traffic flows primarily for model-predictive control purposes and proposing countermeasures

Historically, the motivations behind the development of many continuum models have been a mixed set of criteria from both research directions (Mohammadian, 2021). However, these research directions are not necessarily aligned. That is, a desirable model developed with respect to (a) may, for instance, fall short according to the criteria for a desirable model with respect to (b), and vice versa. This is a point that has been repeatedly revealed in several existing real-world applications of continuum traffic flow models (Kontorinaki et al., 2017; Mohammadian et al., 2021c; Spiliopoulou et al., 2014; Wang et al., 2022). As such, it seems that there might be no ultimate solution to reconcile the criteria for a desirable continuum model for the two research directions in the future.

Considering the challenges and properties of CAVs discussed in Section 4.1, the importance of this point is likely to be more relevant for (mixed) CAV traffic flows. Therefore, we shall consider these two research directions separately, and inspect some inherent challenges for each, and propose recommendations for each case accordingly. Note that our discussions will primarily deal with the real-world traffic flows of commercial CAVs *as they will appear the way they will*, and not the way they should 'ideally' be controlled with respect to certain objectives as in the research direction discussed at the **beginning** of Section 4.

##### 4.2.2.1. *Approach A: Developing models for explanatory purposes.*

In this case, the main focus would be to enable making reproducible inferences about how the control laws of commercial CAVs would affect different physical aspects of traffic flow dynamics as a whole. Some specific questions to consider would be, for instance, how the AI-powered control laws of the commercial CAVs would:

- Alleviate or exacerbate some aspects of traffic flow dynamics regarding congestion emergence, propagation, and dissipation.
- React to small disturbances in various steady-state traffic flows in different settings, and whether CAVs' responses would lead to traffic instabilities or not, and if so to what extent?
- React to the situation in which a single driver in a line of CAVs takes over.

Such questions are also important to be investigated from a safety perspective. Here, we focus on the applications and developments of continuum models for making inferences about (mixed) CAV traffic flow properties, and the criteria that can make models desirable with respect to this end.

To be considered as desirable, it is necessary that the proposed continuum models can:

(1) Replicate the collective interactions of (mixed) CAVs traffic flow similar to the empirical observations to some reasonable extent. For instance, if the focus is to investigate how different headway settings would affect the stability of automated traffic flows, it is imperative that the outcomes from the proposed continuum model matches the expected observations from field testing to some reasonable extent.
(2) Be flexible enough to accommodate at least the major components of the driver-vehicle units (such as flexible headway setting, driver's take-over maneuvers, safety-oriented responses, depending on the focus of inference making).

*4.2.2.1.1. Some challenges for approach A.* One way to investigate whether (1) and (2) are met could be to contrast the Lagrangian variation of the proposed continuum models against the control laws of CAVs, or alternatively, derive the Eulerian continuum model from a derivation basis with sufficient empirical support. In either case, there are two barriers that seem to be hindering the development of such continuum models, as summarized below:

(1) Incompatibilities between CAVs' longitudinal control laws and classical car-following and hyperbolic conservation laws; it seems infeasible to derive a continuum approximation from either hierarchical control laws defined as optimization problems (Guanetti et al., 2018), or AI-powered data-driven rules which govern the longitudinal motions of CAVs (Yu et al., 2021).
(2) Lack of reliable knowledge about the ground truth: CAVs have not been widely introduced to traffic streams yet, and the existing field-testings with CAVs are limited overall and may not satisfy the ground truth criteria with respect to the purpose of continuum model development.

*4.2.2.1.2. Potential solutions to the challenges of approach A.* To address the aforementioned challenges, we propose using a "surrogate" derivation basis drawn from controlled experiments with CAVs. These experiments aim to extract reproducible trajectories, providing empirical support for the proposed models. However, it is important to note that conducting such controlled experiments with human-driven vehicles for the same purpose would be extremely difficult. Unlike the consistent and reproducible results one can obtain with CAVs due to their pre-programmed driving strategies, human-driven vehicles introduce an element of variability that cannot be easily controlled or replicated.

In a broad sense, a surrogate derivation basis with sufficient empirical support could be any alternative microscopic or GKT-based framework that proves useful against certain empirical aspects of (mixed) CAV traffic that are of interest regarding the purpose of model development.

For instance, let us consider continuum model development from a car-following derivation basis. In this case, the equivalence between car-following and continuum models (Jin, 2016), which suggests that for every car-following model presented in terms of an ordinary differential equation (ODE), there exists an equivalent continuum model that captures the collective interactions of identical vehicles precisely the same as the car-following model does. In this scenario, the primary purposes of continuum model development could be to investigate how the variations in certain CAVs' car-following properties (e.g., headway setting, gap-controller modes, collision avoidance) would affect macroscopic traffic flow dynamics (e.g., wave propagation and stability). Such a practice can be enabled by using the procedures given in Jin (2016) to derive continuum approximations of a parsimonious car-following model that meets certain criteria against controlled observations of CAVs' trajectories.

In the following, we discuss how to identify such car-following models using a series of step-by-step systematic evaluations. Note that the primary focus of the discussion is on the microscopic aspects that are likely to have considerable and relevant consequences for collective interactions, and **not** on every minuscule and tiny way in which CAVs differ from human-driven vehicles.

First, comprehensive controlled field testings with CAVs are required in order to extract vehicles' trajectories under various conditions, as pursued in some recent efforts (Li et al., 2022b; Shi and Li, 2021; Stern et al., 2018), and such trajectories can then be considered as a ground truth. It is necessary to design the field testings such that it provides sufficient ground truth with respect to the purpose of model development. For instance, the field testing must result in vehicle trajectories that are 'complete' in the sense of traffic regimes (Sharma et al., 2019a) in order for comprehensive evaluations of the car-following models that are considered as potential candidates to derive the continuum models from (i.e., derivation basis).

Next, a systematic benchmarking analysis is required to investigate whether any of the existing car-following models or their extensions that are transferrable to a continuum framework can adequately replicate the observed CAVs' trajectories, especially the key aspects that are relevant from a continuum modeling perspective. For instance, to design the benchmarking criteria, the factsheets explaining the criteria for a desirable car-following model (Treiber and Kesting, 2013) can be used as a starting point. Note that the focus here is to identify a reasonable derivation basis for continuum approximation, and not the ideal car-following model for modeling individual CAVs, and thus, the investigation focus must be designed accordingly.

For instance, some prominent issues pertaining to car-following model calibration and performance (e.g., overall predictive accuracy, one-to-one trajectory mapping and error distribution (Punzo et al., 2021)) might be of less relevance. Instead, the focus should be given to the capability of the models in capturing the spatio-temporal patterns observed by empirical trajectories as in Treiber and Kesting (2012), in which case, the calibration approach should be adopted accordingly, e.g., by utilizing jam front trajectories, wave speeds, or total time spent as measures rather than conventional car-following performance indicators, such as Root mean squared error (RMSE) of speed, spacing, and acceleration (Punzo et al., 2021; Treiber and Kesting, 2012).

The outcomes of such a systematic evaluation may reveal that to be considered as desirable, the underlying car-following model may have to employ hierarchical car-following control strategies as in recent efforts for modeling CACC driving from empirical observations of CAVs (Milanés and Shladover, 2014, 2016; Shladover et al., 2015). Such car-following models are not quite compatible with the majority of car-following models in the literature that are of continuous ODE-type nature over the speed-spacing domain (Saifuzzaman and Zheng, 2014). For instance, a hierarchical car-following model with multiple control mechanisms might be needed, in which specific mechanisms (e.g., speed control, gap control, safety control.) are activated only within certain stages of car-following (e.g., acceleration, maintaining steady-state, safety control.) as in some recent works (Milanés and Shladover, 2014, 2016; Shladover et al., 2015). Note that in such a car-following model, each mode must still be composed of explicit rules devoid of optimization problems.

Continuum approximation of such a car-following model can then be adopted in making inferences using numerical tests in hypothetical test cases, for example, in order to investigate how real-world CAVs would affect macroscopic traffic flow dynamics.

However, the question about the mathematical tractability and practicality of such continuum models would remain. For example, by incorporating hierarchical driving strategies in the speed adaptation equation, the resulting continuum models will likely be difficult to put into the conservative form (Zhang, 2000), and correspondingly, difficulties may arise in investigating the analytical properties and performing numerical simulations (Abgrall and Karni, 2010). Section 4.2.3 discusses some possible approaches regarding how to overcome such difficulties.

The proposed recommendation above has benefits when the focus is to make inferences about how CAVs' control laws may affect the collective properties of traffic flow. However, the resulting continuum model may not necessarily be well-behaved in terms of calibrations and real-world performance.

*4.2.2.2. Approach B: Developing models for predictive purposes.* Some relevant questions to consider in assessing the proposed model with respect to this aim would be

(1) Can the proposed model provide reasonably accurate estimation of observed real-world (mixed) CAV traffic flows regarding control-oriented criteria (e.g., travel time estimation, delay caused by bottlenecks, maximum length of queue)?
(2) Can the model provide reliable knowledge in terms of the expected changes in traffic flow's operational properties under various scenarios of the penetration rates of CAVs?
(3) Does the model, developed and calibrated on a specific section for a specific time period, perform robustly for another day and traffic pattern?
(4) Is the model well-behaved so that it can provide reliable and practically relevant solutions for all expected scenarios?
(5) Can the model's inherent shortcomings and potential non-physical behaviors be considered insignificant with respect to the purpose of application? And if not, is it feasible to apply ad-hoc treatments to overcome these issues such that the model remains practical for the purpose of application?

These questions can be incorporated quantitively into designing assessment criteria for a given continuum model for (mixed) CAV traffic flows proposed with respect to Approach B. If the model's practical performance for real-world CAV traffic flows is reasonable with respect to such assessment criteria, it could be considered to be desirable for this application purpose, even if the model falls short in some other major aspects of interest in Approach A. Desirable characteristics include: (a) being consistent with the CAVs' underlying control laws, (b) the potential to provide reproducible knowledge about impacts of CAVs' control laws on the macroscopic traffic flow dynamics, from an analytical perspective, and correspondingly (c) the potential to provide reliable knowledge in policy planning and design of CAVs' control laws, (d) the capability to remain mathematically tractable and well-posed under all solution domains, and (e) the inherent capability to provide solutions that are physically relevant for all traffic regimes.

*4.2.2.2.1. Some challenges for approach B.* At first sight, the challenges of this approach may be even more significant compared to Approach A because this problem is two-fold, i.e., how to derive a reasonable continuum model for (mixed) CAV traffic flow to begin with, and how to ensure the proposed model qualifies regarding the assessment criteria discussed above, and these questions could be intertwined and in conflict with one another as discussed at the beginning of this section.

That is, when developing a continuum model from first principles, certain physical considerations (e.g., car-following rules, anisotropy property) are commonly taken into account, which may not have direct implications for ensuring practical relevance and vice versa. Besides, the proportion of studies on real-world applications and assessment of existing continuum models for conventional traffic has been considerably small, where the outcomes are conflicting in certain real-world aspects (e.g., Kontorinaki et al., 2017; Mohammadian et al., 2021c; Spiliopoulou et al., 2014; Wang et al., 2022). Thus, given the lack of sufficient knowledge about the continuum model's real-world performance for the conventional traffic, continuum modeling of real-world (mixed) CAVs for reliable traffic predictions seems to be a far-fetched aim.

One possible solution to enable this approach could be to utilize an existing generic continuum model such as GSOM (Lebacque et al., 2007b) or METANET (Wang et al., 2022) as a base hyperbolic structure, and then further extend it to specify and develop robust mechanisms for the generic functions and terms inside the model since existing real-world applications suggested that the choice of such generic terms would have a considerable impact on the GSOM's real-world performance (Fan and Seibold, 2013; Mohammadian et al., 2021c).

In this case, the main benefit would be to address the question (1) discussed above without being restricted by constraints and glitches of hyperbolic systems that may arise when continuum approximations are to be derived from first principles using car-following or gas-kinetic relations. However, there would still be other important difficulties that might make such a practice seem far-fetched, and for better organization, we summarize these issues into three categories:

(1) How to find a suitable ground as the starting point in order to specify and further extend the generic terms and parameters in the selected continuum model? Some relevant issues to this point would be:
    a. When specifying the generic terms, researchers have conventionally relied on some sorts of plausible assumptions about driving behavior or traffic flow behavior, which may not necessarily result in better model performance.
    b. Adding more parameters and mechanisms might improve numerical accuracy to a certain extent but may not necessarily guarantee model robustness, well-behaved calibration, and reliability for other scenarios.
(2) How to design a model performance evaluation approach with respect to the purpose of model development in this approach in order to enable rigorous assessment of the proposed extended model's performance against the five questions discussed above. The main issues in this case would be:
    a. The objective functions defined for calibration purposes (commonly formulated in terms of squared errors of flow and speed) may not necessarily guarantee the models' best performance for application purposes (e.g., travel time estimations) (Mohammadian et al., 2021c).
    b. The choice of the optimization algorithm can itself affect the model's calibration outcomes and may find the 'optimal' model parameters that would be undesirable with respect to robustness and reliability for other scenarios (Mohammadian et al., 2021c; Spiliopoulou et al., 2017).
(3) How to ensure the measures taken to address (1) and (2) would also guarantee that the model would always remain well-behaved practically and numerically, for all the scenarios of interest to the aim of model development, with respect to the desirability criteria discussed above in terms of five questions. More specifically,
    a. The proposed functional forms and specified mechanisms in the extended model may occasionally yield physically absurd results (e.g., physical instability near jam density, or violation of the anisotropy property), which may still be acceptable

considering the overall model's performance for the scenarios considered.

b. If not treated proactively, such behaviors may dominate the calibration procedure or the model's performance for some application scenarios. At the same time, proposing robust treatments for issues appearing in (a) may not necessarily be feasible and may itself defeat the purpose, and cause issues such as flawed model performance for an otherwise acceptable model behavior.

The main problem is that within conventional frameworks of developing functional forms and model calibration, it might be too difficult to satisfactorily address these three main points with respect to the aim of model development, and as such, it would be difficult to identify the trade-offs. Since the three categories of challenges pertaining to Approach B are intertwined with one another, overcoming them within the classical frameworks of hyperbolic conservation laws, where the all the terms are derived based on physical considerations, may be a far-fetched aim.

*4.2.2.2.2. Potential solutions to the challenge of approach B.* We introduce an alternative hybrid modeling framework, consisting of machine-learning algorithms incorporated into hyperbolic conservation laws, as a potential solution. Hybrid modeling frameworks, consisting of machine-learning and classical statistical models, have been proven useful and recommended in other fields such as traffic safety analysis (Ali et al., 2021; Mannering et al., 2020; Mohammadian et al., 2021a), in order to maximize both the statistical inference and predication capabilities of the proposed models. Recently, combining machine learning and kinematic wave theory has been introduced in some recent works, (Thodi et al., 2022; Yuan et al., 2021), where the focus has been on improving the shortcomings of data-driven AI-powered methods for real-world traffic estimation. Our proposed recommendation would be to adopt this approach, from the other way around, in order to improve the shortcomings of performance of continuum models for real-world traffic flows. Given the complexities and uncertainties surrounding real-world (mixed) CAVs, the feasibility of hybrid AI-powered continuum models could be investigated as a potential solution, and if successful, can bring about numerous benefits such as:

(1) One would not need to face the constraints and glitches of hyperbolic systems that may arise when continuum approximations are to be derived from first principles using car-following or gas-kinetic relations. One could use a predefined generic hyperbolic structure with reasonable performance (such as GSOM) and then utilize machine-learning techniques to further extend and specify the generic mechanisms inside the model such that the model's desirable real-world performance is achieved with respect to the five questions discussed above.
(2) All the theoretical, numerical, and physical constraints for specifying the generic terms (e.g., anisotropy property, finite reaction time, invariant congestion wave speed, numerical stability) can be defined as constraints inside the computational process in the machine-learning algorithms considered.
(3) Due to (2), there would be unprecedented flexibility in identifying the most desirable way to specify and extend the model with respect to real-world traffic. For instance, all the model performance assessment criteria can be defined as constraints in the computational procedure and let the machine-algorithms search for the most desirable way to extend such terms.
(4) Due to (3), there would be increased room for strengthening the calibration and validation approaches with respect to the considerations relevant not only to the model's performance for real-world traffic, but also to the ultimate use of the models' outcomes, i.e., model-predictive control strategies (Ferrara et al., 2018; Han et al., 2021a; Papamichail et al., 2019).

To ensure that the proposed hybrid modeling framework serves the purpose, a systematic step-by-step implementation procedure must be developed, taking into account several considerations.

The first step pertains to the ground truth, i.e., the observations from real-world (mixed) CAV traffic, and their suitability for being considered as benchmarks for evaluating the proposed model. For instance, CAVs will act as probe-vehicles to collect information about traffic conditions along their trajectory, which means that the corresponding data would be Lagrangian in nature. Some researchers have argued that in the case of Lagrangian data, variational formulations of continuum models (i.e., Lagrangian continuum models (Daganzo, 2005; Daganzo, 2006b; Li and Zhang, 2013)) would be of more suitability (Piccoli et al., 2015). While for CAVs, Lagrangian models can certainly be applied straightforwardly, the Lagrangian nature of CAVs' data should not be considered as a challenge for the application of Eulerian continuum models. Note that the reference point is the time–space diagram of traffic states, from which one needs the initial and boundary conditions as the continuum model inputs. Considering this point, we argue that CAVs can provide better opportunities for the application of continuum traffic flow models. This is because the data provided by CAVs can provide a more comprehensive picture about the time–space diagram, which conventionally is reconstructed from loop detectors using smoothing techniques (Mohammadian et al., 2021c; Treiber and Helbing, 2002).

Next, one would have to define the computational procedure for model implementation and calibration, which involves two major components, i.e., (a) the numerical scheme, and (b) optimization algorithm and the model calibration and evaluation design.

Regarding (a), Section 4.2.3. Discusses some research needs that may arise regarding developing new numerical schemes, and here we mainly focus on a few aspects that are more relevant from the model calibration perspective. Some conventional approaches for numerical schemes include approximate Riemann solvers (Kong, 2011; Mohammadian et al., 2021b; Toro, 2013; Zhang et al., 2011), relaxation scheme (Chen et al., 2009; Delis et al., 2014; Mutua, 2016), and discontinuous Galerkin methods (Buli and Xing, 2020; Mohammadian et al., 2021b; Qiao et al., 2017). Since the primary focus is real-world applications, one may use any numerical scheme or its extensions from these categories that works reasonably well for the base hyperbolic model as well as the proposed hybrid model within the scenarios of interest. For real-world implementation, it is necessary that such a numerical scheme is devoid of local spurious oscillations that may arise near strong shockwaves (Mohammadian and Van Wageningen-Kessels, 2018) to avoid potential impacts on the model's calibration performance. However, numerical schemes with increased resolution near shocks (weighted-essentially non-oscillatory (WENO) schemes (Chu et al., 2022; Zhang et al., 2003, 2006; Zhao et al., 2014)) may not necessarily be desirable for real-world implementation, since on the empirical time–space diagram, the observed shockwaves will have a certain width, which is likely to expand if CAVs perform smoother decelerations. In this case, a high-resolution numerical scheme with strong performance near shocks may dominate the model calibration procedure undesirably (Mohammadian et al., 2021b).

Regarding (b), an objective function for calibration purposes and an optimization algorithm should be designed such that the calibration outcomes ensure not only the model's predictive accuracy but also, and perhaps more importantly, the model's robustness, regarding its sensitivity to the changes of boundary conditions, and the trajectory of the jam-front, i.e., the interface between the free-flow and congestion regions. It may turn out that in order to ensure robustness, such aspects might need to be considered directly as objective functions for calibration purposes rather as opposed to the conventional approaches that rely on the squared errors of macroscopic state variables (Kontorinaki et al., 2017; Mohammadian et al., 2021c; Papageorgiou et al., 1990; Porfyri et al., 2016; Spiliopoulou et al., 2014, 2017; Wang et al., 2022). More discussions on this topic can be found elsewhere (Mohammadian et al., 2021c; Ngoduy and Liu, 2007; Ngoduy and Maher, 2012; Treiber and Kesting, 2013 (chapter 16)).

Meanwhile, an objective functions and calibration approach may identify the model's “optimal” parameter set such that the model's *quantitative performance* regarding operational measures (e.g., total time spent), would seem desirable for a sufficiently long time–space domain, but yet, the model's outcomes would be undesirable to be used as inputs for a model-predictive control strategy for shorter time–space intervals (Ferrara et al., 2018; Mohammadian et al., 2021c). For instance, consider a scenario in which compared to the observed data, the calibrated model provides a good estimation of overall traffic states in a time-space window and the corresponding total time spent, which is among the conventional performance measures for control purposes. However, if the congestion region predicted by the model is shifted in the time–space diagram either temporally (e.g., the onset of congestion evolution in the model's outcomes is either considerably sooner or later compared to the observations) or spatially (e.g., the maximum queue predicted by the model is considerably different from the data), such inputs would result in premature or delayed activation of the anticipated control strategy at times of congestion emergence, which may ultimately exacerbate the situation.

#### 4.2.3. Need for utilizing and further development of path-conservative numerical schemes

This section discusses some specific challenges inherent to the numerical approximations of hyperbolic systems, that are likely to be prevalent when CAV-related aspects are adequately incorporated in continuum models.

In Section 4.2.2, we discussed the importance of developing continuum models for (mixed) CAV traffic flows that can reasonably resemble, at the car-following level, the hierarchical driving rules of CAVs, and in such cases, the resulting continuum models are unlikely to be cast into the conservative form as

$$\partial\overrightarrow{U} \Big/ \partial t + \partial F\left(\overrightarrow{U}\right) \Big/ \partial x = S\left(\overrightarrow{U}\right) \tag{4.1}$$

where $\overrightarrow{U} = [U_1, U_2, ..., U_n]$ is the vector of state variables, $F(\overrightarrow{U}) = [F_1, F_2, ..., F_n]$ is the vector of corresponding fluxes, and $S(\overrightarrow{U})$ is the vector of source terms (the terms devoid of gradients). For such a system, the Jacobian matrix is given as $A(\overrightarrow{U}) = \partial F(\overrightarrow{U}) / \partial\overrightarrow{U}$, from which one can obtain the system's characteristic speeds by solving $\det[A(\overrightarrow{U}) - \overrightarrow{\lambda}\overrightarrow{I}] = 0$.

Such models can be presented in the quasi-linear form:

$$\partial\overrightarrow{U} \Big/ \partial t + A\left(\overrightarrow{U}\right)\partial\overrightarrow{U} \Big/ \partial x = S\left(\overrightarrow{U}\right) \tag{4.2}$$

Conventionally, the conservative form in Eq. (4.1) has been considered of utmost importance for studying the model's exact analytic solutions (Lebacque et al., 2007a, 2007b; Mammar et al., 2009; Zhang, 2000, 2001; Zhang et al., 2003). However, the conservative form can only be obtained under certain conditions, such as

(1) When the models are directly derived from the mass conservation principle, e.g., the equilibrium continuum models or phenomenological non-equilibrium models which are directly derived from the flow conservation principle.
(2) When the models are derived from car-following and gas-kinetic theories, and the speed adaptation equation(s) consists of continuous and parsimonious terms that facilitate reformulating the whole system into the conservative form in Eq. (4.1).

However, there have been several existing continuum models for human-driven vehicles that have been derived from more complex behavioral principles that generally cannot be cast into the conservative form, except for some specific cases (e.g., when the Greenshields' fundamental diagram is adopted) (e.g., the models by Berg et al. (2000); Gupta and Katiyar (2006), 2007; Mohammadian et al. (2023); Zhang (1998), 2003. Particularly, some existing continuum models that have desirable physical and behavioral properties are often not globally well-posed for some important and physically relevant scenarios (Helbing, 2009). This issue is likely to be exacerbated if the CAV-related aspects are adequately incorporated into continuum models.

> Remark: Some prospective continuum models with promising properties for (mixed) CAV traffic flows may not be cast into conservative form, but the use of path-conservative approach could help overcoming the inherent numerical difficulties.

A critical challenge caused by the absence of conservative form is that “while for conservative systems, shock relations depend on the solution states to the immediate left/right of the shock, in the non-conservative case they depend not only on those states but also on the path that connects them” (Abgrall and Karni, 2010). The path-conservative theoretical framework, developed in Dal Maso et al. (1995), seems to be the most suitable approach for numerical approximations of non-conservative systems as in Eq. (4.2). However, two main challenges in the path-conservative schemes are: selecting the shock path is somewhat arbitrary, and one cannot investigate whether the obtained solutions are physically relevant (Abgrall and Karni, 2010).

When it comes to continuum traffic flow models, the issues caused by such uncertainties could be alleviated by considering the underlying car-following assumptions when selecting a shock path. The proposed path-conservative schemes can then be assessed against some benchmark cases, where the conservative form is available and solution structures are known (Chu et al., 2022). In recent years, there has been substantial developments in path-conservative numerical schemes (Cao et al., 2022, 2023; Castro Díaz et al., 2019; Dal Maso et al., 1995; Serezhkin and Menshov, 2020), which can also facilitate the future development of continuum models of (mixed) CAVs. Therefore, preserving the conservative property should no longer be considered as a constraint, and the priorities in developing new continuum models of (mixed) CAV traffic should shift towards adequately incorporating the complex human factors and CAVs' driving strategies.

Once a suitable path-conservative scheme has been established for $S(\overrightarrow{U}) = 0$, then it can be generalized to derive well-balanced schemes for $S(\overrightarrow{U}) \neq 0$ in order to approximate the source terms robustly (Cao et al., 2022, 2023; Liu et al., 2015). Unfortunately, numerical studies of continuum traffic flow models, including those for CAVs, have used numerical schemes in which source terms (relaxation terms in particular) have been approximated only based on local observed values and simple time-stepping (Helbing et al., 2001; Helbing and Treiber, 1999; Mohammadian et al., 2021b; Ngoduy, 2012a, 2012b, 2013a, 2013b, 2014; Ngoduy et al., 2004, 2009; Ngoduy and Hoogendoorn, 2003; Ngoduy and Jia, 2017; Treiber et al., 1999; Zheng et al., 2015). Such simple approximations can affect the evolution of disturbances around the steady-state by adding artificial smoothing and diffusions to the numerical results (Liu et al., 2015). This issue is of particular importance for CAV traffic, when highly non-linear continuum models are used and the focus is to make inferences on the system's stability behavior using numerical experiments. Well-balanced path-conservative schemes could facilitate such investigations.

## 5. Conclusions and future directions

Connected and automated vehicles (CAVs) are expected to reshape traffic flow dynamics and present new challenges and opportunities for traffic flow modeling, especially in the foreseeable future when CAVs will be mixed and interacting with human-driven vehicles (Calvert et al., 2018; Mahmassani, 2016). The uncertainties around such mixed traffic flows will also have profound implications in real-world traffic flow modeling, especially the continuum modeling.

We provided a select review of continuum models for conventional human-driven vehicles and revisited some of their properties that likely have implications for CAVs, according to the theoretical and empirical findings as well as the controversies and debates in the literature (Daganzo, 1995, 2006a; Helbing, 2009; Helbing and Johansson, 2009; Helbing et al., 2013; Kontorinaki et al., 2017; Li et al., 2022a; Li and Zhang, 2011; Mohammadian et al., 2021c; Porfyri et al., 2015; Schönhof and Helbing, 2009; Spiliopoulou et al., 2017; Zhang, 1998, 1999, 2000, 2003, 2009; Zhang et al., 2019). Along our review, we also discussed some recent efforts in continuum modeling of (mixed) CAV traffic flow and contrasted them against recent studies on empirical observations and field testings with CAVs (Ciuffo et al., 2021; Gunter et al., 2021; Hu et al., 2022, 2023; Li et al., 2017, 2022b; Milanés and Shladover, 2014; Shi and Li, 2021; Stern et al., 2018; Sun et al., 2020b; Varotto et al., 2015, 2017, 2018; Yu et al., 2021; Zhang et al., 2019).

Overall, it seems the current common practice in the continuum modeling literature is to model CAVs as what they are wished to behave when they are ideally controlled and operated by assuming the maximum utilization of CAVs' inherent resources and capabilities (inter-vehicular communications, platooning, driving at smaller gaps, etc.), in order to improve traffic flow in terms of both efficiency and safety. However, commercial CAVs' driving strategies are likely to be designed conservatively to ensure safety and driver comfort, etc., which are conflicting with traffic efficiency and can lead to undesirable traffic flow properties (Li, 2022), as underscored in recent empirical studies (Ciuffo et al., 2021; Gunter et al., 2021; Hu et al., 2023; Li et al., 2022b). As such, the properties of mixed real-world traffic flows will be much more complex and uncertain.

There is a significant research need to develop new continuum modeling approaches for capturing real-world (mixed) CAV traffic flow properties. The paper highlighted some major problems important for the future development of continuum modeling for real-world (mixed) CAV traffic flows. When developing new continuum models, it is important to differentiate between two major research directions: (a) modeling for explaining purposes, where making reproducible inferences about the physical aspects of macroscopic properties is of the primary interest, and (b) modeling for practical purposes, in which the focus is on the reliable predictions for operation and control. Historically, the motivations behind the development of many continuum models have been a mixed set of criteria from both research directions (Mohammadian, 2021), but it is often the case that many models fall short regarding some major aspects of both directions (Kontorinaki et al., 2017; Mohammadian et al., 2021c; Papageorgiou, 1998), and such shortcomings have sparked controversies in the literature regarding the desirability and consistency of certain modeling frameworks (Daganzo, 1995; Helbing, 2009; Helbing and Johansson, 2009; Schönhof and Helbing, 2009; Zhang, 2009).

The paper highlighted some essential challenges regarding continuum modeling in each research direction and proposed some recommendations to facilitate the future research. As stated by Helbing et al. (2001), "*When aiming at better models, it should be considered that theoretical consistency is just one aspect, and matching empirical data or stylized facts is at least equally important. It is not enough to find partial answers*". We believe this statement is of utmost importance and has profound implications for continuum modeling of real-world (mixed) CAV traffic flow, where numerous considerations related to safety, driver comfort, energy consumption, etc. Will be incorporated into the commercial CAVs' underlying control laws to various degrees (Ciuffo et al., 2021; Li et al., 2017, 2022b; Varotto et al., 2018; Yu et al., 2021).

We hope that this paper can lay the foundation for identifying the trade-offs between continuum models' physical consistencies, mathematical tractability, behavioral plausibility, and performances regarding either their explaining or predicting capabilities, and ultimately act as a catalyst for developing a new paradigm of continuum modeling in the era of CAVs, the transition period in particular.

## Declaration of competing interest

The authors declare that they have no known competing financial interests or personal relationships that could have appeared to influence the work reported in this paper.

## Acknowledgement

The authors are grateful for the numerous insightful and constructive comments received from the anonymous reviewers and the editor, which have substantially improved the quality of this paper. This research was partially funded by the Australian Research Council (ARC) through the Discovery Project (DP210102970) and Dr. Zuduo Zheng's Discovery Early Career Researcher Award (DECRA; DE160100449).

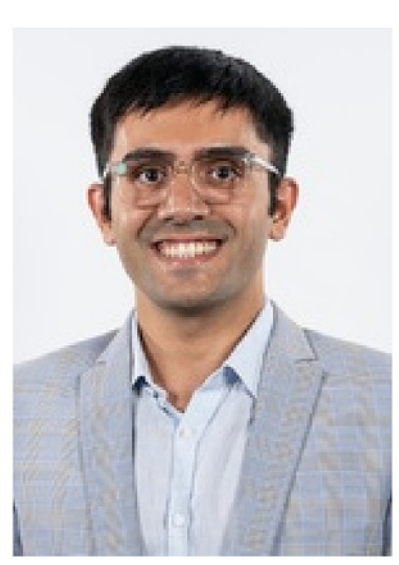

**Saeed Mohammadian** is a Postdoctoral Research Fellow at the School of Civil Engineering, Queensland University of Technology, Australia. He has conducted research in various areas of traffic flow theories, such as development and practical real-world implementation of continuum traffic flow models, numerical methods for continuum models as hyperbolic systems, and safety assessment of macroscopic traffic states using traffic conflict analysis.

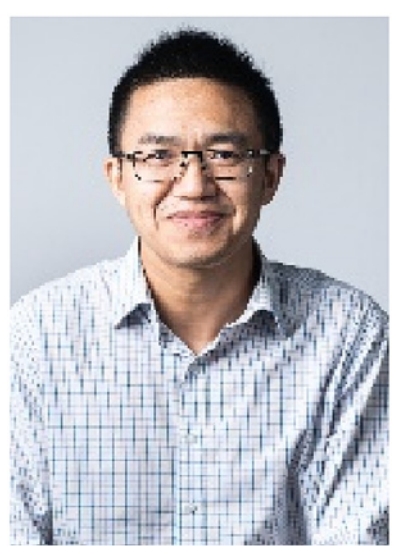

**Zuduo Zheng** is currently TMR TAP Chair (Deputy) and an Associate Professor with the School of Civil Engineering, The University of Queensland. His research interests include traffic flow modeling, travel behavior and decision making, and advanced data analysis techniques, including mathematical modeling, econometrics, numerical optimization in transport engineering, and meta research.

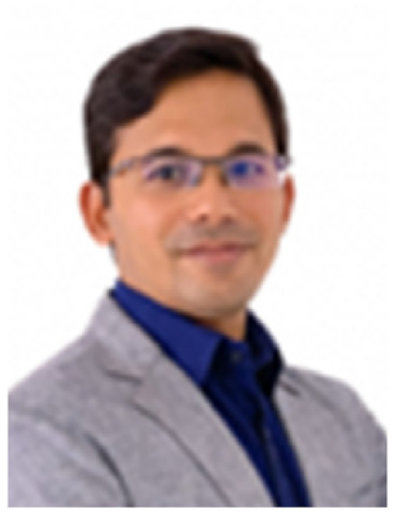

**Md. Mazharul Haque** is currently a Professor at the School of Civil and Environmental Engineering, Queensland University of Technology, Australia. He is a specialist in econometrics and artificial intelligence applications in transport engineering and traffic safety. His research interests include road safety, traffic conflict analysis, connected and automated vehicles, and human factors and driving behaviour.

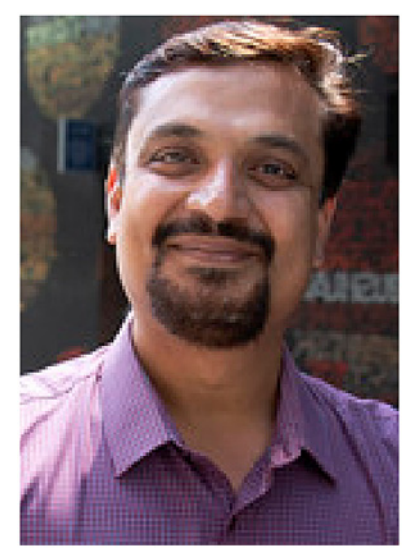

**Ashish Bhaskar** is Professor in Civil Engineering (Transport discipline), at the Queensland University of Technology (QUT), Australia. His expertise includes transport big data analytics, modelling, simulation, and control. He serves as Chief Investigator at QUT Centre for Data science, where he also co-lead Business and Engineering Systems research domain. He holds a Ph.D. degree in Intelligent Transport Systems from Swiss Federal Institute of Technology (EPFL), Lausanne, Switzerland, a M.S. degree in Transport Engineering from the University of Tokyo, Tokyo, Japan, and a B.S. degree of Technology (BTech) degree in Civil Engineering from the Indian Institute of Technology (IIT), Kanpur.